\documentclass[11pt,a4paper,twoside]{amsart}

\RequirePackage{amsthm,amsmath,amsfonts,amssymb}
\RequirePackage[colorlinks,citecolor=blue,urlcolor=blue]{hyperref}

\usepackage{enumitem}
\usepackage{graphics}
\usepackage[mathcal]{euscript}
\usepackage{mathtools}
\usepackage{array}
\usepackage[math]{cellspace}
\usepackage{tikz,tikz-cd}
\usetikzlibrary{fit,shapes.geometric}
\usepackage{xpatch}
\usepackage{thmtools}
\usepackage{thm-restate}
\usepackage{xassoccnt}
\usepackage{ytableau}
\usepackage[font=small,labelfont=bf]{caption}
\usepackage{subfig}
\usepackage{array}
\tikzset{snake it/.style={decorate, decoration={random steps,segment length=8pt,amplitude=2.0pt}}}
\usepackage[capitalise,nameinlink]{cleveref}

\makeatletter
\newcommand*{\centerfloat}{%
  \parindent \z@
  \leftskip \z@ \@plus 1fil \@minus \textwidth
  \rightskip\leftskip
  \parfillskip \z@skip}
\makeatother

\makeatletter
\newcounter{globallistcounter}
\newcounter{thmlistcounter}

\let\realItem\item
\NewDocumentCommand\myItem{o}{%
   \IfNoValueTF{#1}{%
      \realItem%
      \edef\countername{\LastRefSteppedCounter}%
      \def\thegloballistcounter{\csname the\countername\endcsname}%
      \refstepcounter{globallistcounter}%
   }{%
      \realItem[#1]%
      \def\thegloballistcounter{#1}%
      \refstepcounter{globallistcounter}%
   }%
}

\AtBeginEnvironment{enumerate}{%
    \let\item\myItem%
    \stepcounter{globallistcounter}%
}

\def\localItemStringLiteral{LOCALITEM@}
\NewDocumentCommand\myLabel{m}{%
      \realLabel{\localItemStringLiteral#1}%
      \setcounter{thmlistcounter}{\value{globallistcounter}}%
      \def\thethmlistcounter{\theparentcounter\thegloballistcounter}%
      \refstepcounter{thmlistcounter}%
      \realLabel[\parentcountername]{#1}
}

\newlist{thmlist}{enumerate}{1}
\setlist[thmlist]{
    before=\let\item\myItem%
           \let\realLabel\label%
           \let\label\myLabel%
           \stepcounter{globallistcounter}%
           \edef\parentcountername{\LastRefSteppedCounter}%
           \edef\theparentcounter{\csname the\parentcountername\endcsname},
    label=\upshape{(\alph{thmlisti})},
    noitemsep}

\makeatother

\AtEndEnvironment{proof}{\setcounter{claim}{0}}

\declaretheorem[name=Theorem,
                Refname={Theorem,Theorems},
                numberwithin=section,style=plain]{theo}
\declaretheorem[name=Proposition,
               Refname={Proposition,Propositions},
               numberlike=theo,style=plain]{prop}
\declaretheorem[name=Lemma,
               Refname={Lemma,Lemmas},
               numberlike=theo,style=plain]{lemma}

\declaretheorem[name=Definition,
               Refname={Definition,Definitions},
               numberlike=theo,style=plain]{defi}

\declaretheorem[name=Question,
               Refname={Question,Questions},
               numberlike=theo,style=plain]{quest}
\declaretheorem[name=Remark,
               Refname={Remark,Remarkss},
               numberlike=theo,style=remark]{remark}

\newcommand{\bbF}{{\mathbb F}}
\newcommand{\bbZ}{{\mathbb Z}}

\newcommand{\bbR}{{\mathbb R}}
\newcommand{\bbC}{{\mathbb C}}
\newcommand{\calA}{{\mathcal A}}
\newcommand{\calB}{{\mathcal B}}

\newcommand{\bfx}{{\mathbf{x}}}
\newcommand{\bfm}{{\mathbf{m}}}

\newcommand{\bfy}{{\mathbf{y}}}
\newcommand{\bfQ}{{\mathbf{Q}}}

\newcommand{\Aff}{{\rm Aff}}

\newcommand{\GL}{{\rm GL}}

\newcommand{\id}{{\rm id}}

\DeclareMathOperator{\Span}{Span}
\DeclareMathOperator{\im}{im}

\DeclareMathOperator{\Sym}{Sym}
\DeclareMathOperator{\chara}{char}

\newcommand{\trans}[2]{T^{#1}_{#2}{}}
\newcommand{\transgr}[2]{\Gamma^{#1}_{#2}}
\newcommand{\transsimple}[1]{T_{#1}}
\newcommand{\transgrsimple}[1]{\Gamma_{#1}}

\newcommand{\generate}[1]{\langle #1\rangle}
\newcommand{\restrict}[2]{#1|_{#2}}
\newcommand{\ones}{\mathbf{1}}
\newcommand{\demark}[1]{#1\star}
\newcommand{\amat}[4]{\begin{bsmallmatrix}#1 & #2 \\ #3 & #4\end{bsmallmatrix}}
\newcommand{\Amat}[4]{\setlength{\arraycolsep}{2pt}%
\begin{bmatrix*}#1 & #2 \\ #3 & #4\end{bmatrix*}}
\newcommand{\mvec}[2]{\begin{psmallmatrix*}[l]#1 \\ #2\end{psmallmatrix*}}
\newcommand{\Mvec}[2]{\begin{pmatrix*}[l]#1 \\ #2\end{pmatrix*}}
\newcommand{\twoDisjointEdges}{%
\begin{tikzpicture}[scale=0.5]
\draw (0,0)--(1,0);
\draw (2,0)--(3,0);
\filldraw (0,0) circle (3pt);
\filldraw (1,0) circle (3pt);
\filldraw (2,0) circle (3pt);
\filldraw (3,0) circle (3pt);
\end{tikzpicture}%
}

\makeatletter
\@namedef{subjclassname@2020}{\textup{2020} Mathematics Subject Classification}
\makeatother

\begin{document}

\title{Orbits under dual symplectic transvections}
\author{Jonas Sj{\"o}strand}
\address{School of Education, Culture and Communication, Division of Mathematics and Physics, M\"alardalen University, Box~883, 721~23 V\"aster\aa s, Sweden}
\email{jonas.sjostrand@mdu.se}
\date{December, 2023}

\begin{abstract}
Consider an arbitrary field $K$ and
a finite-dimensional vector space $X$ over $K$ equipped with a, possibly degenerate,
symplectic form $\omega$. Given a spanning subset $S$ of $X$,
for each $k$ in $K$ and each vector $s$ in $S$, consider the symplectic transvection
mapping a vector $x$ to $x+k\omega(x,s)s$. The group generated by
these transvections has been extensively studied, and its orbit structure is
known. In this paper, we obtain corresponding results for the orbits of the
dual action on $X^\ast$. As for the non-dual case, the analysis gets harder
when the field contains only two elements. For that field, the dual transvection group is equivalent to a game known as the lit-only sigma game, played on a graph. Our results provide a complete solution to the reachability problem of that game, previously solved only for some special cases.
\end{abstract}

\maketitle

\section{Introduction}\label{sec:introduction}
Let $X$ be a finite-dimensional vector space over a field $K$ and let $\omega$ be
a (possibly degenerate) $K$-valued alternating bilinear form on $X$.
By currying, $\omega$ will also denote the linear mapping from $X$ to the
dual space $X^\ast$ defined by $\omega(x)(y)=\omega(x,y)$.

For any $s\in X$ and any nonzero $k\in K$, let $\transsimple{s,k}$
be the mapping from $X$ to $X$ defined by
\[
\transsimple{s,k}(x) = x + k\omega(x,s)s.
\]
This is called a \emph{(symplectic) transvection}. For notational
convenience we will write $\transsimple{s}$ as a shorthand for $\transsimple{s,1}$
when $K$ has only two elements.

It is easy to see that
$\transsimple{s,k}\circ\transsimple{s,-k}$ is the identity mapping on $X$.
For any subset $S$ of $X$, let the \emph{transvection group}
$\transgrsimple{S}$ be the subgroup of $\GL(X)$ generated by
all $\transsimple{s,k}$ for $s\in S$ and nonzero $k\in K$.
Let $G(S)$ be the (possibly infinite) graph with vertex set $S$ and with an edge between $u$ and $v$ if $\omega(u,v)\ne0$.

The orbit structure of $\transgrsimple{S}$ has been described completely in the literature when $S$ spans $X$ and $G(S)$ is connected. There are three cases that behave differently.

For the case where $K$ has more than two elements, the orbit structure was
found by Brown and Humphries~\cite[Th.~6.5]{BrownHumphriesI1986}.
\begin{theo}[Brown, Humphries 1986]\label{th:brownhumphriesnormalfield}
Suppose $K\ne\bbF_2$. Let $S$ be a spanning subset of $X$
such that $G(S)$ is connected, and consider the group
$\transgrsimple{S}$ acting on $X$. Then, two distinct elements
$x,y\in X$
belong to the same orbit if and only if neither of them belongs to $\ker\omega$.
\end{theo}

When $K=\bbF_2$, the two-element field, we need some more definitions to describe the orbit strucure.
Given a basis $S$ of $X$ and two elements $s,t\in S$ with $\omega(s,t)=1$,
we can construct another basis of $X$ by replacing $t$ with $s+t$.
If a basis $S'$ of $X$ can be obtained from $S$ by a sequence of such
replacements, we say that $S$ and $S'$ are \emph{t-equivalent}, and we also
say that the graphs $G(S)$ and $G(S')$ are {t-equivalent.
A basis $S$ of $X$ is of \emph{orthogonal type} if it
is t-equivalent to a basis $S'$ such that $G(S')$ is a tree that contains
$E_6$ as an induced subgraph. In this case, we also say that the graph $G(S)$ is
of orthogonal type.

If $S$ is a basis of $X$, we denote by $Q_S$ the unique quadratic form
on $X$ such that $Q_S(s)=1$ for any $s\in S$ and
$Q_S(x+y)+Q_S(x)+Q_S(y)=\omega(x,y)$ for any $x,y\in X$.

For a basis of orthogonal type,
Brown and Humphries~\cite[Th.~10.1]{BrownHumphriesII1986} obtained
the following result, independently found by Janssen~\cite{Janssen1983, Janssen1985}.
\begin{theo}[Brown, Humphries 1986; Janssen 1983]
\label{th:brownhumphriesFtwo}
Suppose $K=\bbF_2$, and let $S$ be a
basis of $X$ of orthogonal type.
Then, two elements $x,y\in X\setminus\ker\omega$ belong to the same orbit of $\transgrsimple{S}$ if and only if
$Q_S(x)=Q_S(y)$.
\end{theo}

The orbit structure for a basis not of orthogonal type was not studied
explicitly until twenty years later, when Seven~\cite[Th.~2.6]{Seven2005}
obtained the
following unified description of orbits for any basis.
\begin{theo}[Seven 2005]\label{th:seven}
Suppose $K=\bbF_2$. Let $S$ be a basis of $X$
such that $G(S)$ is connected, and let
$d:X\rightarrow\bbZ_{>0}$ be the function defined as
\[
d(x) = \min\{s : x = x_1 + \dotsb + x_s\ \text{for some}\ x_i\in \transgrsimple{S}S\},
\]
where $\transgrsimple{S}S$ denotes the $\transgrsimple{S}$-orbit containing $S$.
Then, two elements $x,y\in X\setminus\ker\omega$ belong to the same orbit of $\transgrsimple{S}$ if and only if $d(x)=d(y)$.
\end{theo}

Our aim in this paper is to describe the orbit structure of the dual action.
Let $X^\ast$ denote the dual of $X$ as a vector space (forgetting
about the bilinear form) and, for each $g\in\GL(X)$, define
the dual mapping $g^\ast\in\GL(X^\ast)$ by letting
$g^\ast(\alpha) = \alpha\circ g$.
The duals of the elements of $\transgrsimple{S}$ form a subgroup of $\GL(X^\ast)$
denoted by $\transgrsimple{S}^\ast$, a \emph{dual transvection group}.
\begin{quest}\label{qu:main}
For a spanning subset $S$ of $X$,
when do $\alpha,\beta\in X^\ast$ belong to the same orbit of $\transgrsimple{S}^\ast$?
\end{quest}
The problem splits into the same
three cases as the non-dual version, and we will find dual analogues to each one
of the three theorems above. 

The paper is organized as follows. First, in \cref{sec:sigmagame},
we review previous work and discuss an alternative description of a dual
transvection group
$\transgrsimple{S}^\ast$ when $K=\bbF_2$ as a game played on a graph. In \cref{sec:results}, we present our main results. In 
\cref{sec:prerequisites} we introduce some tools and notation.
Then, in \cref{sec:components,sec:lineardependence}
we show that \cref{qu:main} can be reduced to the case where $G(S)$ is connected
and, if $S$ is finite, to the case where $S$ is a basis of $X$. In \cref{sec:generalfield}, we derive some group isomorphism results that will be
essential for our analysis, both when $K\ne\bbF_2$ and when $K=\bbF_2$.
After that, we are finally ready to prove our main results, and it is done in three
sections: \cref{sec:morethantwo} for the case where $K\ne\bbF_2$, \cref{sec:orthogonal} for bases of orthogonal type and \cref{sec:linegraphs}
for other bases in the $K=\bbF_2$ case. Finally, in \cref{sec:nondual}, we fill a gap in the theory of ordinary (non-dual) transvection groups by a theorem
about the case where $G(S)$ is not connected.

Our ambition is to make the presentation as self-contained as possible.

\section{Previous work}\label{sec:sigmagame}
To the best of our knowledge, no one has studied the orbit structure
of dual transvection groups over arbitrary fields. However, for the
$K=\bbF_2$ case and with $S$ being a basis of $X$,
there are many related results in the literature, using
varying terminology.
For instance, Janssen~\cite{Janssen1983,Janssen1985} refers to
$\transgrsimple{S}$ as a \emph{monodromy group},
to its orbits (together with $X$ and $\omega$) as \emph{vanishing lattices}
and to $S$ as a \emph{weakly distinguished basis}.

For the case where $K=\bbF_2$ and $S$ is a basis
of orthogonal type, Shapiro et.~al.~\cite[Lemma~4.6]{ShapiroEtAl1998} found the
number of $\transgrsimple{S}^\ast$-orbits, but they did not address
\cref{qu:main}.

Many authors have studied dual transvection groups over $\bbF_2$ in terms of a
one-player game called the \emph{lit-only sigma game}.
It is played
on an undirected graph, each vertex of which has a lamp that is either
on or off. A move consists of choosing any lit vertex, that is, a vertex
whose lamp is on, and toggle the state of all adjacent vertices. Usually, the
goal is to reach a position with as few lit vertices as possible. (Note, that
it is impossible to turn off all lamps since a move always leaves the played vertex unaffected.)

The graph in our case is $G(S)$, where $S$ is a basis of $X$,
and the game state is an element
$\alpha\in X^\ast$, where a vertex $v\in S$ is lit if and only $\alpha(v)=1$.
Playing a lit vertex $s\in S$ corresponds to the dual transvection $\transsimple{s}^\ast$ acting on $\alpha$ to reach the new state $\alpha\circ\transsimple{s}$.
Clearly,
\[
(\alpha\circ\transsimple{s})(v) = \alpha(v + \omega(v,s)s)
= \begin{cases}
\alpha(v) + 1 & \text{if $v$ is adjacent to $s$,} \\
\alpha(v) & \text{otherwise,}
\end{cases}
\]
so the rules of the lit-only sigma game are followed.
Also, for a non-lit vertex $s\in S$, the transvection $\transsimple{s}$ leaves $\alpha$
unaffected.

Conversely, given a simple graph $G=(V,E)$ on which to play the lit-only sigma game, we can choose $X$ as the vector space freely generated by $V$ and define
$\omega$ by letting $\omega(u,v)=1$ for $u,v\in V$ if and only if $(u,v)\in E$.
Then, $G(V)$ is isomorphic to $G$, so we have converted the lit-only sigma game
to the dual transvection group $\transgrsimple{V}^\ast$.

The lit-only sigma game is also equivalent to Mozes's \emph{game of numbers} \cite{Mozes1990} played with coefficients in $\bbF_2$ rather than in $\bbR$
or $\bbC$. This is a linear representation of the simply-laced Coxeter group given by the graph. We refer to \cite[Ch.~4]{BjornerBrentiBook2005} for the details.

\Cref{qu:main} becomes equivalent to the reachability problem for the game:
Given two game positions, can one position be reached from the other by a sequence of moves?

The lit-only sigma game was originally obtained from adding the lit-only rule
to a game called a \emph{$\sigma^-$-automaton}, introduced by Sutner~\cite{Sutner1989}.
It has been studied extensively, with a focus on the
minimum and maximum number of lit vertices that can be obtained
\cite{WangWu2007, GoldwasserKlostermeyer2009, GoldwasserWangWu2009,
WangWu2009, WangWu2010, GoldwasserWangWu2011, Huang2013}.
In 2008, Huang and Weng \cite{HuangWeng2008} solved the reachability problem for
the lit-only sigma game for graphs of type A, D and E in the classification of irreducible Coxeter groups. Soon thereafter,
Wu~\cite{Wu2009} and Huang and Weng~\cite{HuangWeng2010}
studied the game on line graphs of simple graphs, and our discussion in \cref{sec:linegraphs} is based heavily on their approach.
In 2015, Huang \cite{Huang2015} solved the reachability problem for graphs
whose corresponding alternating form $\omega$ is nondegenerate; a characterization
of these graphs was given by Reeder \cite{Reeder2005}.
Finally, in 2020, Vorstermans \cite{Vorstermans2020} studied the group structure
of the lit-only sigma game and a generalization of it.

After finishing the research for the present paper,
we found a nicely written introductory text by Wu and Xiang \cite{WuXiang2019} mentioning
results similar to ours for the case $K=\bbF_2$ but without any mathematical argument. For the proofs, the authors referred to a paper that is ``to appear'' ([52] in their bibliography list), but it does not seem to have appeared yet.

\section{Main results}\label{sec:results}
In \cref{sec:components}, we will show that \cref{qu:main} can
be reduced to the case where $G(S)$ is connected, so that will be an assumption
for the rest of this section.

Our first main result is an analogue to \cref{th:brownhumphriesnormalfield}.
\begin{restatable*}{theo}{maindualnormalfield}\label{th:maindualnormalfield}
Suppose $K\ne\bbF_2$, and let $S$ be a spanning subset of $X$
such that $G(S)$ is connected. Then, two nonzero elements
$\alpha,\beta\in X^\ast$
belong to the same orbit of $\transgrsimple{S}^\ast$ if and only if
$\beta-\alpha\in\im\omega$.
\end{restatable*}

In \cref{sec:lineardependence}, we will show that, if $S$ is finite, \cref{qu:main} can be reduced to the case where $S$ is a basis of $X$.
Our second main result deals with bases of orthogonal type and
is an analogue to \cref{th:brownhumphriesFtwo}.
\begin{restatable*}{theo}{maindualmodtwoorthogonal}\label{th:maindualmodtwoorthogonal}
Suppose $K=\bbF_2$, and
let $S$ be a basis of $X$ of orthogonal type.
Then, two nonzero elements $\alpha$ and $\beta$ of $X^\ast$
belong to the same orbit of $\transgrsimple{S}^\ast$ if and only if
there is an $x\in X$ such that $\omega(x)=\alpha+\beta$ and
$Q_S(x)=\alpha(x)$.
\end{restatable*}

\begin{remark}
A priori, checking whether there is an $x\in X$ such that $\omega(x)=\alpha+\beta$ and $Q_S(x)=\alpha(x)$ might be computationally hard. In fact, this turns out to be an easy task:
Finding an $x\in X$ such that $\omega(x)=\alpha+\beta$ is just a matter of solving a linear system of equations, and if there is such an $x$, we do not have to compute
the $(Q_S+\alpha)$-value for all such $x$. By \cref{lm:qorbits}, we only need
to check whether there is an $x_0\in\ker\omega$ with $(Q_S+\alpha)(x_0)=1$,
which is easy since $Q_S+\alpha$ is a linear function when restricted to
$\ker\omega$. If no such $x_0$ exists, all $x$ in $\omega^{-1}(\alpha+\beta)$
have the same $(Q_S+\alpha)$-value, so we only need to compute it once.
\end{remark}

Our third and final main result concerns bases not of orthogonal type.
It uses the following beautiful theorem by Cuypers \cite[Th.~3.3 and Th.~3.4]{Cuypers2021}.
\begin{theo}[Cuypers 2021]\label{th:cuypers}
A connected graph $G$ is the line graph of some connected multigraph if and and only
if $G$ is not of orthogonal type.
\end{theo}
Here, by a \emph{multigraph} we mean a graph where multiple edges are allowed but
not loops, and by the \emph{line graph} of a multigraph $(V,E)$ we mean the simple
graph whose vertex set is $E$ and where there is an edge between $e_1$ and $e_2$
if they have exactly one endpoint in common.

Suppose $K=\bbF_2$.
Let $S$ be a basis of $X$ not of orthogonal type and suppose that $G(S)$ is connected. Then, by \cref{th:cuypers}, the graph $G(S)$ is the
line graph of some connected multigraph $G=(V,E)$, so we can identify $S$ with $E$.

Let $\generate{V}$ be the free vector space over $K$ on the set $V$ and equip
$\generate{V}$ with an inner product $\omega_{\generate{V}}(\cdot,\cdot)$ such that
$V$ is an orthogonal basis.
Define a linear mapping $\partial$, called the \emph{boundary map} from $X$ to
$\generate{V}$ by letting
$\partial$ of an edge be the sum of its endpoints, and let the adjoint mapping
$\delta$, called the \emph{co-boundary map},
from $\generate{V}$ to $X^\ast$ be defined
by $\delta(y)(x) = \omega_{\generate{V}}(y,\partial(x))$.

For an element $y$ in $\generate{V}$, let $d_0(y)$ and $d_1(y)$ denote the number of
zero and one coordinates of $y$ in the basis $V$, respectively,
and let $d(y)=\min\{d_0(y),d_1(y)\}$.

\begin{restatable*}{theo}{maindualmodtwolinegraph}\label{th:maindualmodtwolinegraph}
Suppose $K=\bbF_2$. Let $S$ be a basis of $X$ not of orthogonal type, and suppose $G(S)$ is connected.
Then two elements $\beta,\gamma\in X^\ast$ belong to the same orbit of
$\transgrsimple{S}^\ast$ if and only if $\gamma-\beta\in\im\omega$
and either $\beta\not\in\im\delta$ or $\beta\in\im\delta$ and
$d(y)=d(z)$ for some (or, equivalently, any) $y,z\in\generate{V}$ such that
$\delta(y)=\beta$ and $\delta(z)=\gamma$.
\end{restatable*}

\begin{remark}
By another result of Cuypers
\cite[Th.~1.1]{Cuypers2021}, a connected (ordinary) graph is the line graph of a
multigraph if and only if it does not contain any of the 32 graphs in \cref{fig:thirtytwographs} as an induced subgraph.
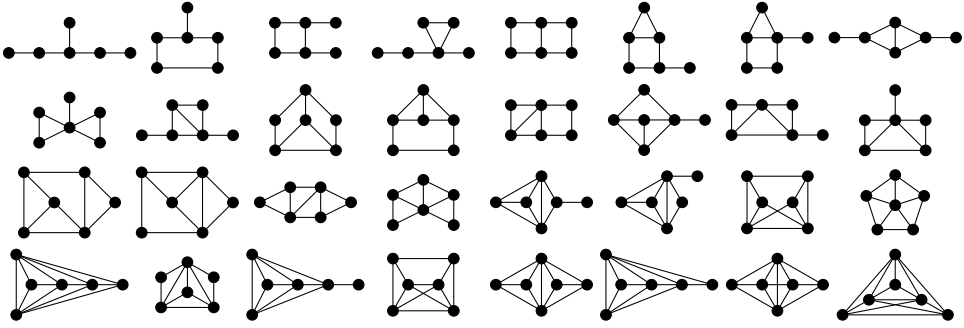
\begin{figure}
\def\scale{0.04}
\def\radius{50pt}
\def\colwidth{12mm}
\centerfloat
\begin{tabular}{%
>{\centerfloat\arraybackslash}m{\colwidth}%
>{\centerfloat\arraybackslash}m{\colwidth}%
>{\centerfloat\arraybackslash}m{\colwidth}%
>{\centerfloat\arraybackslash}m{\colwidth}%
>{\centerfloat\arraybackslash}m{\colwidth}%
>{\centerfloat\arraybackslash}m{\colwidth}%
>{\centerfloat\arraybackslash}m{\colwidth}%
>{\centerfloat\arraybackslash}m{\colwidth}%
}
\begin{tikzpicture}[scale=\scale]
\filldraw (-20,  0) coordinate (A) circle (\radius);
\filldraw (-10,  0) coordinate (B) circle (\radius);
\filldraw (  0,  0) coordinate (C) circle (\radius);
\filldraw ( 10,  0) coordinate (D) circle (\radius);
\filldraw ( 20,  0) coordinate (E) circle (\radius);
\filldraw (  0, 10) coordinate (F) circle (\radius);
\path[-] (A) edge (B);
\path[-] (B) edge (C);
\path[-] (C) edge (D) edge (F);
\path[-] (D) edge (E);
\end{tikzpicture}
&
\begin{tikzpicture}[scale=\scale]
\filldraw (-10,  0) coordinate (A) circle (\radius);
\filldraw (-10, 10) coordinate (B) circle (\radius);
\filldraw (  0, 10) coordinate (C) circle (\radius);
\filldraw ( 10, 10) coordinate (D) circle (\radius);
\filldraw ( 10,  0) coordinate (E) circle (\radius);
\filldraw (  0, 20) coordinate (F) circle (\radius);
\path[-] (A) edge (B) edge (E);
\path[-] (B) edge (C);
\path[-] (C) edge (D) edge (F);
\path[-] (D) edge (E);
\end{tikzpicture}
&
\begin{tikzpicture}[scale=\scale]
\filldraw ( 20,  0) coordinate (A) circle (\radius);
\filldraw ( 10,  0) coordinate (B) circle (\radius);
\filldraw (  0,  0) coordinate (C) circle (\radius);
\filldraw (  0, 10) coordinate (D) circle (\radius);
\filldraw ( 10, 10) coordinate (E) circle (\radius);
\filldraw ( 20, 10) coordinate (F) circle (\radius);
\path[-] (A) edge (B);
\path[-] (B) edge (C) edge (E);
\path[-] (C) edge (D);
\path[-] (D) edge (E);
\path[-] (E) edge (F);
\end{tikzpicture}
&
\begin{tikzpicture}[scale=\scale]
\filldraw (-20,  0) coordinate (A) circle (\radius);
\filldraw (-10,  0) coordinate (B) circle (\radius);
\filldraw (  0,  0) coordinate (C) circle (\radius);
\filldraw ( 10,  0) coordinate (D) circle (\radius);
\filldraw ( -5, 10) coordinate (E) circle (\radius);
\filldraw (  5, 10) coordinate (F) circle (\radius);
\path[-] (A) edge (B);
\path[-] (B) edge (C);
\path[-] (C) edge (D) edge (E) edge (F);
\path[-] (E) edge (F);
\end{tikzpicture}
&
\begin{tikzpicture}[scale=\scale]
\filldraw ( 20,  0) coordinate (A) circle (\radius);
\filldraw ( 10,  0) coordinate (B) circle (\radius);
\filldraw (  0,  0) coordinate (C) circle (\radius);
\filldraw (  0, 10) coordinate (D) circle (\radius);
\filldraw ( 10, 10) coordinate (E) circle (\radius);
\filldraw ( 20, 10) coordinate (F) circle (\radius);
\path[-] (A) edge (B) edge (F);
\path[-] (B) edge (C) edge (E);
\path[-] (C) edge (D);
\path[-] (D) edge (E);
\path[-] (E) edge (F);
\end{tikzpicture}
&
\begin{tikzpicture}[scale=\scale]
\filldraw ( 20,  0) coordinate (A) circle (\radius);
\filldraw ( 10,  0) coordinate (B) circle (\radius);
\filldraw (  0,  0) coordinate (C) circle (\radius);
\filldraw (  0, 10) coordinate (D) circle (\radius);
\filldraw (  5, 20) coordinate (E) circle (\radius);
\filldraw ( 10, 10) coordinate (F) circle (\radius);
\path[-] (A) edge (B);
\path[-] (B) edge (C) edge (F);
\path[-] (C) edge (D);
\path[-] (D) edge (E) edge (F);
\path[-] (E) edge (F);
\end{tikzpicture}
&
\begin{tikzpicture}[scale=\scale]
\filldraw ( 20, 10) coordinate (A) circle (\radius);
\filldraw ( 10,  0) coordinate (B) circle (\radius);
\filldraw (  0,  0) coordinate (C) circle (\radius);
\filldraw (  0, 10) coordinate (D) circle (\radius);
\filldraw (  5, 20) coordinate (E) circle (\radius);
\filldraw ( 10, 10) coordinate (F) circle (\radius);
\path[-] (A) edge (F);
\path[-] (B) edge (C) edge (F);
\path[-] (C) edge (D);
\path[-] (D) edge (E) edge (F);
\path[-] (E) edge (F);
\end{tikzpicture}
&
\begin{tikzpicture}[scale=\scale]
\filldraw (-20,  0) coordinate (A) circle (\radius);
\filldraw (-10,  0) coordinate (B) circle (\radius);
\filldraw (  0,  5) coordinate (C) circle (\radius);
\filldraw (  0, -5) coordinate (D) circle (\radius);
\filldraw ( 10,  0) coordinate (E) circle (\radius);
\filldraw ( 20,  0) coordinate (F) circle (\radius);
\path[-] (A) edge (B);
\path[-] (B) edge (C) edge (D);
\path[-] (C) edge (D) edge (E);
\path[-] (D) edge (E);
\path[-] (E) edge (F);
\end{tikzpicture}
\\
\begin{tikzpicture}[scale=\scale]
\filldraw (-10,  5) coordinate (A) circle (\radius);
\filldraw (-10, -5) coordinate (B) circle (\radius);
\filldraw (  0,  0) coordinate (C) circle (\radius);
\filldraw ( 10,  5) coordinate (D) circle (\radius);
\filldraw ( 10, -5) coordinate (E) circle (\radius);
\filldraw (  0, 10) coordinate (F) circle (\radius);
\path[-] (A) edge (B) edge (C);
\path[-] (B) edge (C);
\path[-] (C) edge (D) edge (E) edge (F);
\path[-] (D) edge (E);
\end{tikzpicture}
&
\begin{tikzpicture}[scale=\scale]
\filldraw (-10,  0) coordinate (A) circle (\radius);
\filldraw (  0,  0) coordinate (B) circle (\radius);
\filldraw (  0, 10) coordinate (C) circle (\radius);
\filldraw ( 10, 10) coordinate (D) circle (\radius);
\filldraw ( 10,  0) coordinate (E) circle (\radius);
\filldraw ( 20,  0) coordinate (F) circle (\radius);
\path[-] (A) edge (B);
\path[-] (B) edge (C) edge (E);
\path[-] (C) edge (D) edge (E);
\path[-] (D) edge (E);
\path[-] (E) edge (F);
\end{tikzpicture}
&
\begin{tikzpicture}[scale=\scale]
\filldraw (-10,  0) coordinate (A) circle (\radius);
\filldraw (-10, 10) coordinate (B) circle (\radius);
\filldraw (  0, 20) coordinate (C) circle (\radius);
\filldraw ( 10, 10) coordinate (D) circle (\radius);
\filldraw ( 10,  0) coordinate (E) circle (\radius);
\filldraw (  0, 10) coordinate (F) circle (\radius);
\path[-] (A) edge (B) edge (E) edge (F);
\path[-] (B) edge (C);
\path[-] (C) edge (D) edge (F);
\path[-] (D) edge (E);
\path[-] (E) edge (F);
\end{tikzpicture}
&
\begin{tikzpicture}[scale=\scale]
\filldraw (-10,  0) coordinate (A) circle (\radius);
\filldraw (-10, 10) coordinate (B) circle (\radius);
\filldraw (  0, 20) coordinate (C) circle (\radius);
\filldraw ( 10, 10) coordinate (D) circle (\radius);
\filldraw ( 10,  0) coordinate (E) circle (\radius);
\filldraw (  0, 10) coordinate (F) circle (\radius);
\path[-] (A) edge (B) edge (E);
\path[-] (B) edge (C) edge (F);
\path[-] (C) edge (D) edge (F);
\path[-] (D) edge (E) edge (F);
\end{tikzpicture}
&
\begin{tikzpicture}[scale=\scale]
\filldraw ( 20,  0) coordinate (A) circle (\radius);
\filldraw ( 10,  0) coordinate (B) circle (\radius);
\filldraw (  0,  0) coordinate (C) circle (\radius);
\filldraw (  0, 10) coordinate (D) circle (\radius);
\filldraw ( 10, 10) coordinate (E) circle (\radius);
\filldraw ( 20, 10) coordinate (F) circle (\radius);
\path[-] (A) edge (B) edge (F);
\path[-] (B) edge (C) edge (E);
\path[-] (C) edge (D) edge (E);
\path[-] (D) edge (E);
\path[-] (E) edge (F);
\end{tikzpicture}
&
\begin{tikzpicture}[scale=\scale]
\filldraw ( 20,  0) coordinate (A) circle (\radius);
\filldraw ( 10,  0) coordinate (B) circle (\radius);
\filldraw (  0,-10) coordinate (C) circle (\radius);
\filldraw (-10,  0) coordinate (D) circle (\radius);
\filldraw (  0, 10) coordinate (E) circle (\radius);
\filldraw (  0,  0) coordinate (F) circle (\radius);
\path[-] (A) edge (B);
\path[-] (B) edge (C) edge (E) edge (F);
\path[-] (C) edge (D) edge (F);
\path[-] (D) edge (E) edge (F);
\end{tikzpicture}%
&
\begin{tikzpicture}[scale=\scale]
\filldraw (-10,  5) coordinate (A) circle (\radius);
\filldraw (-10, -5) coordinate (B) circle (\radius);
\filldraw (  0,  5) coordinate (C) circle (\radius);
\filldraw ( 10,  5) coordinate (D) circle (\radius);
\filldraw ( 10, -5) coordinate (E) circle (\radius);
\filldraw ( 20, -5) coordinate (F) circle (\radius);
\path[-] (A) edge (B) edge (C);
\path[-] (B) edge (C) edge (E);
\path[-] (C) edge (D) edge (E);
\path[-] (D) edge (E);
\path[-] (E) edge (F);
\end{tikzpicture}
&
\begin{tikzpicture}[scale=\scale]
\filldraw (-10,  5) coordinate (A) circle (\radius);
\filldraw (-10, -5) coordinate (B) circle (\radius);
\filldraw (  0,  5) coordinate (C) circle (\radius);
\filldraw ( 10,  5) coordinate (D) circle (\radius);
\filldraw ( 10, -5) coordinate (E) circle (\radius);
\filldraw (  0, 15) coordinate (F) circle (\radius);
\path[-] (A) edge (B) edge (C);
\path[-] (B) edge (C) edge (E);
\path[-] (C) edge (D) edge (E) edge (F);
\path[-] (D) edge (E);
\end{tikzpicture}
\\
\begin{tikzpicture}[scale=\scale]
\filldraw (  0,  0) coordinate (A) circle (\radius);
\filldraw (  0, 20) coordinate (B) circle (\radius);
\filldraw ( 20, 20) coordinate (C) circle (\radius);
\filldraw ( 30, 10) coordinate (D) circle (\radius);
\filldraw ( 20,  0) coordinate (E) circle (\radius);
\filldraw ( 10, 10) coordinate (F) circle (\radius);
\path[-] (A) edge (B) edge (E) edge (F);
\path[-] (B) edge (C) edge (F);
\path[-] (C) edge (D) edge (E);
\path[-] (D) edge (E);
\path[-] (E) edge (F);
\end{tikzpicture}
&
\begin{tikzpicture}[scale=\scale]
\filldraw (  0,  0) coordinate (A) circle (\radius);
\filldraw (  0, 20) coordinate (B) circle (\radius);
\filldraw ( 20, 20) coordinate (C) circle (\radius);
\filldraw ( 30, 10) coordinate (D) circle (\radius);
\filldraw ( 20,  0) coordinate (E) circle (\radius);
\filldraw ( 10, 10) coordinate (F) circle (\radius);
\path[-] (A) edge (B) edge (E);
\path[-] (B) edge (C) edge (F);
\path[-] (C) edge (D) edge (E) edge (F);
\path[-] (D) edge (E);
\path[-] (E) edge (F);
\end{tikzpicture}
&
\begin{tikzpicture}[scale=\scale]
\filldraw (  0,  0) coordinate (A) circle (\radius);
\filldraw ( 10,  5) coordinate (B) circle (\radius);
\filldraw ( 20,  5) coordinate (C) circle (\radius);
\filldraw ( 30,  0) coordinate (D) circle (\radius);
\filldraw ( 20, -5) coordinate (E) circle (\radius);
\filldraw ( 10, -5) coordinate (F) circle (\radius);
\path[-] (A) edge (B) edge (F);
\path[-] (B) edge (C) edge (F);
\path[-] (C) edge (D) edge (E) edge (F);
\path[-] (D) edge (E);
\path[-] (E) edge (F);
\end{tikzpicture}
&
\begin{tikzpicture}[scale=\scale]
\filldraw (-10,  5) coordinate (A) circle (\radius);
\filldraw (-10, -5) coordinate (B) circle (\radius);
\filldraw (  0,  0) coordinate (C) circle (\radius);
\filldraw ( 10,  5) coordinate (D) circle (\radius);
\filldraw ( 10, -5) coordinate (E) circle (\radius);
\filldraw (  0, 10) coordinate (F) circle (\radius);
\path[-] (A) edge (B) edge (C) edge (F);
\path[-] (B) edge (C);
\path[-] (C) edge (D) edge (E) edge (F);
\path[-] (D) edge (E) edge (F);
\end{tikzpicture}
&
\begin{tikzpicture}[scale=\scale]
\filldraw ( 15,  0) coordinate (A) circle (\radius);
\filldraw (  5,  0) coordinate (B) circle (\radius);
\filldraw (  0, {-10*sqrt(3)/2}) coordinate (C) circle (\radius);
\filldraw (  0, {10*sqrt(3)/2}) coordinate (D) circle (\radius);
\filldraw (-15,  0) coordinate (E) circle (\radius);
\filldraw ( -5,  0) coordinate (F) circle (\radius);
\path[-] (A) edge (B);
\path[-] (B) edge (C) edge (D);
\path[-] (C) edge (D) edge (E) edge (F);
\path[-] (D) edge (E) edge (F);
\path[-] (E) edge (F);
\end{tikzpicture}
&
\begin{tikzpicture}[scale=\scale]
\filldraw ( 10,{10*sqrt(3)/2}) coordinate (A) circle (\radius);
\filldraw (  5,  0) coordinate (B) circle (\radius);
\filldraw (  0, {-10*sqrt(3)/2}) coordinate (C) circle (\radius);
\filldraw (  0, {10*sqrt(3)/2}) coordinate (D) circle (\radius);
\filldraw (-15,  0) coordinate (E) circle (\radius);
\filldraw ( -5,  0) coordinate (F) circle (\radius);
\path[-] (A) edge (D);
\path[-] (B) edge (C) edge (D);
\path[-] (C) edge (D) edge (E) edge (F);
\path[-] (D) edge (E) edge (F);
\path[-] (E) edge (F);
\end{tikzpicture}
&
\begin{tikzpicture}[scale=\scale]
\filldraw ( -5,  0) coordinate (A) circle (\radius);
\filldraw (-10,-{10*sqrt(3)/2}) coordinate (B) circle (\radius);
\filldraw (-10, {10*sqrt(3)/2}) coordinate (C) circle (\radius);
\filldraw (  5,  0) coordinate (D) circle (\radius);
\filldraw ( 10,-{10*sqrt(3)/2}) coordinate (E) circle (\radius);
\filldraw ( 10, {10*sqrt(3)/2}) coordinate (F) circle (\radius);
\path[-] (A) edge (B) edge (C) edge (E);
\path[-] (B) edge (C) edge (D) edge (E);
\path[-] (C) edge (F);
\path[-] (D) edge (E) edge (F);
\path[-] (E) edge (F);
\end{tikzpicture}
&
\begin{tikzpicture}[scale=\scale]
\filldraw (  0,  0) coordinate (A) circle (\radius);
\filldraw (  0, 10) coordinate (B) circle (\radius);
\filldraw ({10*sin(360/5)},  {10*cos(360/5)}) coordinate (C) circle (\radius);
\filldraw ({10*sin(2*360/5)},{10*cos(2*360/5)}) coordinate (D) circle (\radius);
\filldraw ({10*sin(3*360/5)},{10*cos(3*360/5)}) coordinate (E) circle (\radius);
\filldraw ({10*sin(4*360/5)},{10*cos(4*360/5)}) coordinate (F) circle (\radius);
\path[-] (A) edge (B) edge (C) edge (D) edge (E) edge (F);
\path[-] (B) edge (C) edge (F);
\path[-] (C) edge (D);
\path[-] (D) edge (E);
\path[-] (E) edge (F);
\end{tikzpicture}
\\
\begin{tikzpicture}[scale=\scale]
\filldraw ( 40,  0) coordinate (A) circle (\radius);
\filldraw ( 30,  0) coordinate (B) circle (\radius);
\filldraw ( 20,  0) coordinate (C) circle (\radius);
\filldraw ( 10,  0) coordinate (D) circle (\radius);
\filldraw (  5,-10) coordinate (E) circle (\radius);
\filldraw (  5, 10) coordinate (F) circle (\radius);
\path[-] (A) edge (B) edge (E) edge (F);
\path[-] (B) edge (C) edge (E) edge (F);
\path[-] (C) edge (D) edge (E) edge (F);
\path[-] (D) edge (E) edge (F);
\path[-] (E) edge (F);
\end{tikzpicture}
&
\begin{tikzpicture}[scale=\scale]
\filldraw (  0,  0) coordinate (A) circle (\radius);
\filldraw (-{5*sqrt(3)}, -5) coordinate (B) circle (\radius);
\filldraw (-{5*sqrt(3)},  5) coordinate (C) circle (\radius);
\filldraw ( {5*sqrt(3)}, -5) coordinate (D) circle (\radius);
\filldraw ( {5*sqrt(3)},  5) coordinate (E) circle (\radius);
\filldraw (  0, 10) coordinate (F) circle (\radius);
\path[-] (A) edge (B) edge (D) edge (F);
\path[-] (B) edge (C) edge (D) edge (F);
\path[-] (C) edge (F);
\path[-] (D) edge (E) edge (F);
\path[-] (E) edge (F);
\end{tikzpicture}
&
\begin{tikzpicture}[scale=\scale]
\filldraw ( 40,  0) coordinate (A) circle (\radius);
\filldraw ( 30,  0) coordinate (B) circle (\radius);
\filldraw ( 20,  0) coordinate (C) circle (\radius);
\filldraw ( 10,  0) coordinate (D) circle (\radius);
\filldraw (  5,-10) coordinate (E) circle (\radius);
\filldraw (  5, 10) coordinate (F) circle (\radius);
\path[-] (A) edge (B);
\path[-] (B) edge (C) edge (E) edge (F);
\path[-] (C) edge (D) edge (E) edge (F);
\path[-] (D) edge (E) edge (F);
\path[-] (E) edge (F);
\end{tikzpicture}
&
\begin{tikzpicture}[scale=\scale]
\filldraw ( -5,  0) coordinate (A) circle (\radius);
\filldraw (-10,-{10*sqrt(3)/2}) coordinate (B) circle (\radius);
\filldraw (-10, {10*sqrt(3)/2}) coordinate (C) circle (\radius);
\filldraw (  5,  0) coordinate (D) circle (\radius);
\filldraw ( 10,-{10*sqrt(3)/2}) coordinate (E) circle (\radius);
\filldraw ( 10, {10*sqrt(3)/2}) coordinate (F) circle (\radius);
\path[-] (A) edge (B) edge (C) edge (D) edge (E);
\path[-] (B) edge (C) edge (D) edge (E);
\path[-] (C) edge (F);
\path[-] (D) edge (E) edge (F);
\path[-] (E) edge (F);
\end{tikzpicture}
&
\begin{tikzpicture}[scale=\scale]
\filldraw ( 15,  0) coordinate (A) circle (\radius);
\filldraw (  5,  0) coordinate (B) circle (\radius);
\filldraw (  0, {-10*sqrt(3)/2}) coordinate (C) circle (\radius);
\filldraw (  0, {10*sqrt(3)/2}) coordinate (D) circle (\radius);
\filldraw (-15,  0) coordinate (E) circle (\radius);
\filldraw ( -5,  0) coordinate (F) circle (\radius);
\path[-] (A) edge (B) edge (C) edge (D);
\path[-] (B) edge (C) edge (D);
\path[-] (C) edge (D) edge (E) edge (F);
\path[-] (D) edge (E) edge (F);
\path[-] (E) edge (F);
\end{tikzpicture}
&
\begin{tikzpicture}[scale=\scale]
\filldraw ( 40,  0) coordinate (A) circle (\radius);
\filldraw ( 30,  0) coordinate (B) circle (\radius);
\filldraw ( 20,  0) coordinate (C) circle (\radius);
\filldraw ( 10,  0) coordinate (D) circle (\radius);
\filldraw (  5,-10) coordinate (E) circle (\radius);
\filldraw (  5, 10) coordinate (F) circle (\radius);
\path[-] (A) edge (B) edge (F);
\path[-] (B) edge (C) edge (E) edge (F);
\path[-] (C) edge (D) edge (E) edge (F);
\path[-] (D) edge (E) edge (F);
\path[-] (E) edge (F);
\end{tikzpicture}
&
\begin{tikzpicture}[scale=\scale]
\filldraw ( 15,  0) coordinate (A) circle (\radius);
\filldraw (  5,  0) coordinate (B) circle (\radius);
\filldraw (  0, {-10*sqrt(3)/2}) coordinate (C) circle (\radius);
\filldraw (  0, {10*sqrt(3)/2}) coordinate (D) circle (\radius);
\filldraw (-15,  0) coordinate (E) circle (\radius);
\filldraw ( -5,  0) coordinate (F) circle (\radius);
\path[-] (A) edge (B) edge (C) edge (D);
\path[-] (B) edge (C) edge (D) edge (F);
\path[-] (C) edge (D) edge (E) edge (F);
\path[-] (D) edge (E) edge (F);
\path[-] (E) edge (F);
\end{tikzpicture}
&
\begin{tikzpicture}[scale=\scale]
\filldraw (  0,  0) coordinate (A) circle (\radius);
\filldraw (-{5*sqrt(3)}, -5) coordinate (B) circle (\radius);
\filldraw (-{10*sqrt(3)},-10) coordinate (C) circle (\radius);
\filldraw ( {5*sqrt(3)}, -5) coordinate (D) circle (\radius);
\filldraw ({10*sqrt(3)},-10) coordinate (E) circle (\radius);
\filldraw (  0, 10) coordinate (F) circle (\radius);
\path[-] (A) edge (B) edge (D) edge (F);
\path[-] (B) edge (C) edge (D) edge (E) edge (F);
\path[-] (C) edge (D) edge (E) edge (F);
\path[-] (D) edge (E) edge (F);
\path[-] (E) edge (F);
\end{tikzpicture}
\end{tabular}
\caption{The 32 graphs t-equivalent to $E_6$.}\label{fig:thirtytwographs}
\end{figure}
(Cuypers accidentally included
the graph
\smash{%
\parbox{\wd0}{\hbox{%
\def\scale{0.03}%
\def\radius{50pt}%
\begin{tikzpicture}[scale=\scale]
\filldraw ( 20,  0) coordinate (A) circle (\radius);
\filldraw ( 10,  0) coordinate (B) circle (\radius);
\filldraw (  0, -8) coordinate (C) circle (\radius);
\filldraw (-10,  0) coordinate (D) circle (\radius);
\filldraw (  0,  8) coordinate (E) circle (\radius);
\filldraw (  0,  0) coordinate (F) circle (\radius);
\path[-] (A) edge (B);
\path[-] (B) edge (C) edge (E) edge (F);
\path[-] (C) edge (D) edge (F);
\path[-] (D) edge (E);
\path[-] (E) edge (F);
\end{tikzpicture}
}}}
in \cite[$E_6^\mathrm{(15)}$ in Fig.~1]{Cuypers2021} which is the line graph of
\parbox{\wd0}{\hbox{%
\def\scale{0.03}%
\def\radius{50pt}%
\begin{tikzpicture}[scale=\scale]
\filldraw (  0,  0) coordinate (A) circle (\radius);
\filldraw ( 10,  0) coordinate (B) circle (\radius);
\filldraw ( 20,  0) coordinate (C) circle (\radius);
\filldraw ( 30,  5) coordinate (D) circle (\radius);
\filldraw ( 30, -5) coordinate (E) circle (\radius);
\filldraw ( 40, -5) coordinate (F) circle (\radius);
\path[-] (A) edge (B);
\path[-] (B) edge[bend left=35] (C);
\path[-] (B) edge[bend right=35] (C);
\path[-] (C) edge (D) edge (E);
\path[-] (E) edge (F);
\end{tikzpicture}%
}}%
.)
To use \cref{th:maindualmodtwolinegraph}, we must also find the multigraph $G$ of which $G(S)$ is a line graph. Cuypers presents an efficient algorithm for that in
\cite[Sec.~3]{Cuypers2021Whitney}.
\end{remark}

\section{Tools and notation}\label{sec:prerequisites}
In this section we introduce some notation and terminology and recall some basic results that will be used throughout the paper. There are three subsections: one about permutation groups, one about affine transvections and one about quadratic forms.

\subsection{Permutation groups}\label{sec:permutationgroups}
Given a (possibly infinite) set $X$, a \emph{permutation group on $X$} is a subgroup of $\Sym(X)$.
In other words, it is a group whose elements are
permutations of $X$ and whose multiplication is
function composition.

The set of fixed points under the action of a permutation group $G$ on $X$
is denoted by $X^G$.

Given a permutation group $G$ on $X$ and a $G$-invariant subset
$Y\subseteq X$,
the restriction to $Y$ of the elements of $G$ form a permutation group
on $Y$ denoted by $\restrict{G}{Y}$. More generally, if $G$ acts on a set $Z$ by
$G\xrightarrow{\rho}\Sym(Z)$, where $\rho$ is a group homomorphism,
we let $\restrict{G}{Z}$ denote the permutation group $\im\rho$ on $Z$.

For convenience, if $\psi$ is a map and $F$ is a set of maps,
we introduce the notation $\psi F$ for $\{\psi\circ f:f\in F\}$
and $F\psi$ for $\{f\circ\psi:f\in F\}$.
\begin{defi}
Let $G$ and $H$ be permutation groups on $X$ and $Y$, respectively,
and let $\psi$ be a bijective mapping from $X$ to $Y$ such that $\psi G=H\psi$.
Then, the group isomorphism $\xi$ from $G$ to $H$ given by $\xi(g)=\psi\circ g\circ\psi^{-1}$ is said to be \emph{induced by $\psi$}, and we say that
$G$ and $H$ are \emph{isomorphic as permutation groups}.
\end{defi}

To check that a mapping induces a group isomorphism, it suffices to check
that it commutes with the action of group generators, as the following lemma
entails.
\begin{lemma}\label{lm:generatorssuffice}
Let $G$ and $H$ be permutation groups on $X$ and $Y$, respectively,
and let $A$ and $B$ be generating subsets of $G$ and $H$, respectively,
both closed under inversion.
Let $\psi$ be a mapping from $X$ to $Y$ such that $\psi A=B\psi$.
Then $G$ acts on the quotient set $X/\psi:=\{\psi^{-1}(y):y\in\im\psi\}$ of nonempty fibers,
$\im\psi$ is an $H$-invariant subset of $Y$, and
$\psi$ (seen as a map from $X/\psi$ to $\im\psi$) induces a group isomorphism from
$\restrict{G}{X/\psi}$ to $\restrict{H}{\im\psi}$.
\end{lemma}
\begin{proof}
Let us first show that $\psi G=H\psi$.

Take any $g\in G$ and write $g=g_n\circ\dotsb\circ g_1$ for some 
$g_1,\dotsc,g_n\in A$. Since $\psi A=B\psi$, there are $h_1,\dotsc,h_n\in B$ such
that every square in the following diagram commutes.
\[
\begin{tikzcd}
    X \arrow{r}{g_1} \arrow{d}{\psi}
    & X \arrow{r}{g_2} \arrow{d}{\psi}
    & \cdots \arrow{r}{g_n} \arrow{d}{\psi}
    & X \arrow{d}{\psi} \\
    Y \arrow{r}{h_1}
    & Y \arrow{r}{h_2}
    & \cdots \arrow{r}{h_n}
    & Y
\end{tikzcd}
\]
Thus, $\psi\circ g = h_n\circ\dotsb\circ h_1\circ\psi\in H\psi$, so
$\psi G\subseteq H\psi$. Similarly, take any $h\in H$ and write $h=h_n\circ\dotsb\circ h_1$ for some $h_1,\dotsc,h_n\in B$. Since $\psi A=B\psi$, there are
$g_1,\dotsc,g_n\in A$ such that the diagram above commutes.
Thus, $h\circ\psi = \psi\circ g_n\circ\dotsb\circ g_1\in \psi G$, so
$\psi G\supseteq H\psi$. We conclude that $\psi G=H\psi$.

From $\psi G=H\psi$ it follows immediately that $\im\psi$ is $H$-invariant.

Take any $g\in G$ and any fiber $p=\psi^{-1}(y)$, where $y\in Y$.
Since $\psi G=H\psi$, there is an $h\in H$ such that
$(\psi\circ g)(\psi^{-1}(y)) = (h\circ\psi)(\psi^{-1}(y)) = \{h(y)\}$,
so $g(\psi^{-1}(y))\subseteq\psi^{-1}(h(y))$. Thus, for any fiber $p$,
any $g\in G$ maps the elements of $p$ into a single fiber $q$.
But $g^{-1}$ maps $q$ into a single fiber too, namely $p$, so
$g$ takes fibers to fibers, and $G$ acts on the quotient set $X/\psi$.

Since $\psi$ is bijective as a map from $X/\psi$ to $\im\psi$, it follows
that $\psi$ induces a group isomorphism from
$\restrict{G}{X/\psi}$ to $\restrict{H}{\im\psi}$.
\end{proof}

\subsection{Affine transvections}\label{sec:notation}
Recall that a bilinear form $\omega$ on a vector space $X$ is
\emph{skew-symmetric} if
$\omega(x,y)=-\omega(y,x)$ for any $x,y\in X$ and \emph{alternating} if
$\omega(x,x)=0$ for any $x\in X$. If $\chara K\ne2$, these
notions coincide, but in characteristic two we need to distinguish between them.

Let $X$ be a finite-dimensional vector space over $K$ equipped with a
skew-symmetric bilinear
form $\omega$.
In our analysis we will need to extend the notion of transvections and
transvection groups as follows.

For any $a\in K$, nonzero $k\in K$ and $s\in X$,
we define the \emph{affine transvection} $\trans{a}{s,k}$ to be the mapping
from $X$ to $X$ defined by
\[
\trans{a}{s,k}(x) = x + k(\omega(x,s)+a)s.
\]
For notational convenience we write $\trans{a}{s}$ as a shorthand for $\trans{a}{s,1}$ when $K=\bbF_2$.

Note that if $\omega(s,s)=0$ then
\begin{equation}\label{eq:inverse}
\begin{aligned}
& (\trans{a}{s,-k}\circ\trans{a}{s,k})(x)= \\
&=x + k(\omega(x,s)+a)s - k[\omega(x + k[\omega(x,s)+a]s,s)+a]s \\
&=x + k(\omega(x,s)+a)s - k(\omega(x,s)+a)s \\
&=x
\end{aligned}
\end{equation}

Now, let $S$ be a set, let $\alpha$ be a mapping from $S$ to $K$ and let $\phi$ be a mapping from $S$ to $X$ such that $\omega(\phi(s),\phi(s))=0$ for any $s\in S$.
Then, we define the \emph{affine transvection group} $\transgr{\alpha}{S,\phi}$
to be the subgroup of $\Aff(X)$ generated by all $\trans{\alpha(s)}{\phi(s),k}$
for $s\in S$ and $k\in K\setminus\{0\}$. (This is a group by virtue of \cref{eq:inverse}.)

For notational convenience, if $S\subseteq X$ and $\phi$ is the identity map on $S$, we omit
$\phi$ and write $\transgr{\alpha}{S}$ as a shorthand for $\transgr{\alpha}{S,\id_S}$.

\subsection{Quadratic forms}
We need to recall some theory about quadratic forms over a field with only two elements.

Let $X$ be a finite-dimensional vector space over the two-element field $\bbF_2$.
A \emph{quadratic form} $Q$ on $X$ is a mapping from $X$ to $\bbF_2$ such that
$Q(x+y)+Q(x)+Q(y)$ is a bilinear function of $x$ and $y$. Given an ordered
basis $S=\{s_1,\dotsc,s_n\}$ of $X$, each element $x\in X$ can be
written as a column vector $\bfx$ with the $S$-coordinates of $x$.
In this basis, each quadratic form $Q$ corresponds to an lower-triangular
$n$-by-$n$ matrix $\bfQ$
such that $Q(x)=\bfx^T\bfQ\bfx$.
The bilinear form $\omega(x,y)=Q(x+y)+Q(x)+Q(y)$
corresponds to the skew-symmetric matrix $\bfQ+\bfQ^T$ (with zeros on the diagonal)
since $\omega(x,y)=\bfx^T(\bfQ+\bfQ^T)\bfy$. The off-diagonal elements
of $\bfQ$ can be recovered from $\bfQ+\bfQ^T$, so from $\bfQ+\bfQ^T$
together with the diagonal elements $\bfQ_{i,i}=Q(s_i)$ it is possible
to recover $Q$.
In particular, for any
alternating bilinear form $\omega$ on $X$ and any basis $S$ of $X$, there
is a unique quadratic form $Q_S$ such that $Q_S(s)=1$ for any $s\in S$ and
$Q_S(x+y)+Q_S(x)+Q_S(y)=\omega(x,y)$ for any $x,y\in X$.

There is a combinatorial interpretation of $Q(x)$ in terms of the
graph $G(S)$, namely that $Q(x)$ is, modulo two, the number of vertices plus the
number of edges in the subgraph of $G(S)$ induced by the vertices that sum to $x$.
In other words, it is the \emph{Euler characteristic} modulo two of this induced subgraph.

\section{Handling multiple components}\label{sec:components}
In this section we show that \cref{qu:main} can be reduced to the case where
$G(S)$ is connected. A corresponding theorem for the non-dual case is given in
\cref{sec:nondual}.

As usual, let $X$ be a finite-dimensional vector space over $K$ equipped with an alternating bilinear form $\omega$.
\begin{theo}\label{th:components}
Let $S$ be a spanning subset of $X$, let $\{S_i\}_{i\in I}$ be the connected components of $G(S)$ and let $X_i=\Span(S_i)$. Then the following holds.
\begin{itemize}
\item
For each $i\in I$, the restriction map $\restrict{\cdot}{X_i}$ gives a group
homomorphism from $\transgrsimple{S}$ to $\restrict{\transgrsimple{S_i}}{X_i}$, and the family of these maps is an isomorphism
$\transgrsimple{S}\cong \prod_{i\in I}\restrict{\transgrsimple{S_i}}{X_i}$.
\item
Two $\alpha,\beta\in X^\ast$ belong to the same orbit of
$\transgrsimple{S}^\ast$ if and only if $\restrict{\alpha}{X_i}$ and
$\restrict{\beta}{X_i}$ belong to the same orbit of
$(\restrict{\transgrsimple{S_i}}{X_i})^\ast$ for any $i\in I$.
\end{itemize}
\end{theo}
\begin{proof}
$X_i$ is $\transsimple{s,k}$-invariant
for any $s\in S_i$ and also for any $s\in S\setminus S_i$ since then
$\omega(x,s)=0$ for any $x\in X_i$.
It follows that $X_i$ is $\transgrsimple{S}$-invariant and that
the restriction map
$\restrict{\cdot}{X_i}$ is a group homomorphism
from $\transgrsimple{S}$ to $\restrict{\transgrsimple{S_i}}{X_i}$.

To show that $\transgrsimple{S}\cong\prod_{i\in I}\restrict{\transgrsimple{S_i}}{X_i}$,
let $P$ be any group, and let
$\{\phi_i:P\rightarrow\restrict{\transgrsimple{S_i}}{X_i}\}_{i\in I}$ be a
family of group homomorphisms. We need to show that there is a unique
homomorphism $\phi:P\rightarrow\transgrsimple{S}$ such that
$\restrict{\phi(p)}{X_i} = \phi_i(p)$ for any $p\in P$ and any $i\in I$.
Since $S$ spans $X$ we can choose a basis $B\subseteq S$ of $X$.
To define $\phi(p)$ it is enough to specify it on $B$.

In order to satisfy
$\restrict{\phi(p)}{X_i} = \phi_i(p)$ for any $p\in P$ and any $i\in I$,
we have to define $\phi(p)$ such that
$\phi(p)(b) := \phi_i(p)(b)$, where $i$ is the unique element
in $I$ such that $b\in S_i$. To check that this single possible candidate is good
enough, consider any $x=\sum_{b\in B}\lambda_b b\in X_i$. Then $x=x_1+x_2$ where $x_1=\sum_{b\in B\cap S_i}\lambda_b b$
and $x_2=\sum_{b\in B\setminus S_i}\lambda_b b$.
Since $x$ and $x_1$ belong to $X_i$, so does $x_2$, and it follows that
$\omega(x_2,s)=0$ for any $s\in S\setminus S_i$ and thus for any $s\in S$.
This implies that $\phi_i(p)(x_2)=\phi(p)(x_2)=x_2$ for any $s\in X_i$,
so
\begin{multline*}
\phi_i(p)(x)=\phi_i(p)(x_1+x_2)=x_2+\phi_i(p)(x_1)
=x_2 + \sum_{b\in B\cap S_i}\lambda_b \phi_i(p)(b)\\
=x_2 + \sum_{b\in B\cap S_i}\lambda_b \phi(p)(b)
=\phi(p)(x).
\end{multline*}
We conclude that
$\restrict{\phi(p)}{X_i} = \phi_i(p)$ for any $p\in P$ and any $i\in I$.

Now let $\alpha,\beta\in X^\ast$. Suppose there is a $g\in\transgrsimple{S}$
such that $\alpha\circ g=\beta$. Then
$(\restrict{\alpha}{X_i})\circ(\restrict{g}{X_i}) = \restrict{\beta}{X_i}$ for any $i\in I$. Conversely, suppose instead that, for any $i\in I$,
there is a $g_i\in\restrict{\transgrsimple{S_i}}{X_i}$ such that $\restrict{\alpha}{X_i}\circ g_i = \restrict{\beta}{X_i}$. Then, by the direct product result,
there is a (unique) $g\in\transgrsimple{S}$ such that $\restrict{g}{X_i}=g_i$ for each $i\in I$.
Hence, $(\restrict{\alpha}{X_i})\circ(\restrict{g}{X_i}) = \restrict{\beta}{X_i}$, so $\alpha\circ g$ coincides with
$\beta$ on $X_i$ for any $i\in I$. Since $S=\bigcup_{i\in I}S_i$
spans $X$, we conclude that $\alpha\circ g=\beta$.
\end{proof}

\begin{remark}
We note that
Brown and Humphries~\cite[Prop.~2.8]{BrownHumphriesI1986} showed that
$\transgrsimple{S}$ is the direct product of the subgroups $\transgrsimple{S_i}$,
without restricting to $X_i$.
\end{remark}

\section{Handling linear dependence}\label{sec:lineardependence}
In this section we show that, if $S$ is finite,
\cref{qu:main} can be reduced to the case where
$S$ is a basis of $X$. The approach is essentially a dual variant of what
Brown and Humphries call ``extensions of symplectic spaces''
in Section~6 of \cite{BrownHumphriesI1986}.

Let $X$ be a finite-dimensional vector space over $K$ equipped
with an alternating bilinear form $\omega$.
Let $S$ be a finite subset of $X$,
and let $Y$ be the free vector space over $K$
on a set $B=\{b_s\}_{s\in S}$ of symbols indexed by $S$. Equip
$Y$ with an alternating bilinear form $\omega_Y$
defined by $\omega_Y(b_s,b_t)=\omega(s,t)$
for any $s,t\in S$.
Let $p$ be the linear map from $Y$ to $X$ defined by $p(b_s)=s$
for any $s\in S$. This map clearly preserves the bilinear form,
that is, $\omega(p(y_1),p(y_2)) = \omega_Y(y_1,y_2)$ for any
$y_1,y_2\in Y$. The dual map $p^\ast$ from $X^\ast$ to $Y^\ast$
is defined by $p^\ast(\alpha)=\alpha\circ p$.
\begin{theo}
Suppose $S$ is finite and spans $X$. Then, $\im p^\ast$ is $\transgrsimple{B}^\ast$-invariant and $p^\ast$ induces a group
isomorphism between
$\transgrsimple{S}^\ast$ and $\restrict{\transgrsimple{B}^\ast}{\im p^\ast}$
In particular,
$\alpha,\beta\in X^\ast$ belong to the same orbit of $\transgrsimple{S}^\ast$ if and only if
$p^\ast(\alpha)$ and $p^\ast(\beta)$ belong to the same orbit
of $\transgrsimple{B}^\ast$.
\end{theo}
\begin{proof}
Since $S$ spans $X$, the map $p$ is surjective.
We claim that the diagram
\[
\begin{tikzcd}
    Y \arrow{r}{\transsimple{b_s,k}} \arrow[two heads]{d}{p} & Y \arrow[two heads]{d}{p} \\
    X \arrow{r}{\transsimple{s,k}} & X
\end{tikzcd}
\]
commutes for any $s\in S$ and any nonzero $k\in K$.
Indeed, since $p$ preserves the bilinear form, for any $y\in Y$ we have
\begin{multline*}
p(\transsimple{b_s,k}(y))
= p(y + k\omega_Y(y,b_s)b_s)
= p(y) + k\omega(p(y),p(b_s))p(b_s)\\
= p(y) + k\omega(p(y),s)s
= \transsimple{s,k}(p(y)).
\end{multline*}
This implies that the dual diagram
\[
\begin{tikzcd}
    Y^\ast & Y^\ast \arrow{l}{\transsimple{b_s,k}^\ast} \\
    X^\ast \arrow[hook]{u}{p^\ast} & X^\ast \arrow[hook]{u}{p^\ast} \arrow{l}{\transsimple{s,k}^\ast}
\end{tikzcd}
\]
commutes as well. By \cref{lm:generatorssuffice}, $\im p^\ast$ is
$\transgrsimple{B}^\ast$-invariant and $p^\ast$ induces a group isomorphism from
$\restrict{\transgrsimple{S}^\ast}{X^\ast/p^\ast}$ to
$\restrict{\transgrsimple{B}^\ast}{\im p^\ast}$. Finally, since
$p$ is surjective, $p^\ast$ is injective, so
$\restrict{\transgrsimple{S}^\ast}{X^\ast/p^\ast} = \transgrsimple{S}^\ast$.
\end{proof}

\section{Main lemmas}\label{sec:generalfield}
In this section, we derive some group isomorphism results that
will be essential for our analysis, both when $K\ne\bbF_2$ and when $K=\bbF_2$.
The first lemma relates a dual transvection group to an affine one.
\begin{lemma}\label{lm:pressingstolampsmaster}
Let $X$ and $Y$ be finite-dimensional vector spaces over $K$
equipped with skew-symmetric bilinear forms
$\omega_X$ and $\omega_Y$, and let $\phi$ be a linear mapping from
$X$ to $Y$ respecting the bilinear forms, that is, for any $x_1,x_2\in X$
it holds that $\omega_Y(\phi(x_1),\phi(x_2))=\omega_X(x_1,x_2)$.
Let $\theta$ be the linear mapping from $Y$ to $X^\ast$ defined by
$\theta(y)=\omega_Y(y)\circ\phi$, and,
for any $\alpha\in X^\ast$,
let $\theta_\alpha$ be the affine mapping from $Y$ to $X^\ast$ defined by
$\theta_\alpha(y) = \theta(y) + \alpha$.

Then, for any subset $S\subseteq X$ such that $\omega_X(s,s)=0$ for any $s\in S$,
$\transgr{\alpha}{S,\phi}$ acts on $Y/\ker\theta$, $\im\theta_\alpha$ is
$\transgrsimple{S}^\ast$-invariant and
$\theta_\alpha$ induces a group isomorphism from
$\restrict{\transgr{\alpha}{S,\phi}}{Y/\ker\theta}$ to $\restrict{\transgrsimple{S}^\ast}{\im\theta_\alpha}$.
\end{lemma}
\begin{proof}
By \cref{lm:generatorssuffice}, it suffices to show that the diagram
\[
\begin{tikzcd}
    Y \arrow{r}{\trans{\alpha(s)}{\phi(s),k}} \arrow{d}{\theta_\alpha} & Y \arrow{d}{\theta_\alpha} \\
    X^\ast \arrow{r}{\transsimple{s,-k}^\ast} & X^\ast
\end{tikzcd}
\]
commutes for any $s\in S$ and any nonzero $k\in K$.
This is a straightforward matter of verification:
\begin{multline*}
\transsimple{s,-k}^\ast(\theta_\alpha(y))(x)
= \theta_\alpha(y)(\transsimple{s,-k}(x))
= \bigl((\omega_Y(y)\circ\phi) + \alpha\bigr)\bigl(x - k\omega_X(x,s)s\bigr)
\\= \bigl((\omega_Y(y)\circ\phi) + k[\omega_Y(y,\phi(s))+\alpha(s)]\omega_X(s) + \alpha\bigr)(x)
\\= \theta_\alpha(y + k[\omega_Y(y,\phi(s))+\alpha(s)]\phi(s))(x)
= \theta_\alpha(\trans{\alpha(s)}{\phi(s),k}(y))(x).
\end{multline*}
The second last equality follows from
$\omega_Y(\phi(s))\circ\phi = \omega_X(s)$, which in turn follows from
$\omega_Y(\phi(s),\phi(x)) = \omega_X(s,x)$.
\end{proof}

The following specialization of \cref{lm:pressingstolampsmaster}
will come to use in \cref{sec:orthogonal}.
\begin{lemma}\label{lm:pressingstolamps}
Let $\theta_\alpha$ be the affine mapping from $X$ to $X^\ast$ defined by
$\theta_\alpha(x) = \omega(x) + \alpha$. Then, for any subset $S\subseteq X$,
it holds that $\transgr{\alpha}{S}$ acts on $X/\ker\omega$, $\im\theta_\alpha$
is $\transgrsimple{S}^\ast$-invariant, and
$\theta_\alpha$ induces a group isomorphism from
$\restrict{\transgr{\alpha}{S}}{X/\ker\omega}$ to $\restrict{\transgrsimple{S}^\ast}{\im\theta_\alpha}$.
\end{lemma}
\begin{proof}
This follows from \cref{lm:pressingstolampsmaster} with $Y=X$ and
$\phi$ the identity mapping on $X$.
\end{proof}

The third lemma is much more special but will be a key ingredient in induction proofs later on.
\begin{lemma}\label{lm:groupactionisomorphism}
Let $X$ be a finite-dimensional vector space over $K$ equipped with a skew-symmetric bilinear form $\omega$,
and let $Y$ be a subspace of $X$.
Let $x$ be an element of $X$ and
let $\psi$ be the affine transformation on $X$ that adds $x$ to its argument. 

Then, for any $\alpha\in X^\ast$ and any subset
$S$ of $Y$ such that $\omega(s,s)=0$ for any $s\in S$, it holds that
$\psi$ induces a group isomorphism between
$\restrict{\transgr{\alpha+\omega(x)}{S}}{Y}$ and
$\restrict{\transgr{\alpha}{S}}{x+Y}$.
\end{lemma}
\begin{proof}
By \cref{lm:generatorssuffice}, it suffices to show that
the following diagram commutes for any
$s\in S$, any $a\in K$ and any nonzero $k\in K$.
\[
\begin{tikzcd}
    Y \arrow{r}{\trans{a+\omega(x,s)}{s,k}} \arrow{d}{\psi} & Y  \arrow{d}{\psi} \\
    x+Y \arrow{r}{\trans{a}{s,k}} & x+Y
\end{tikzcd}
\]
This is a straightforward matter of verification:
\begin{multline*}
\psi(\trans{a+\omega(x,s)}{s,k}(y)) = \trans{a+\omega(x,s)}{s,k}(y)+x = y + x + k(\omega(y+x,s) + a)s
\\= \trans{a}{s,k}(y+x)
= \trans{a}{s,k}(\psi(y)).
\end{multline*}
\end{proof}

\section{The case \texorpdfstring{$K\ne\bbF_2$}{the field has more than two elements}}\label{sec:morethantwo}
The goal of this section is to prove our first main result, \cref{th:maindualnormalfield}.

To handle the case where $S$ is infinite, we need a couple of general lemmas
about infinite graphs. These are certainly well known, but we could not find a
proper reference.
\begin{lemma}\label{lm:infiniteconnectivity}
In a connected (possibly infinite) graph $G=(V,E)$, for any finite subset $S$ of $V$
there is a connected finite induced subgraph of $G$ whose vertex set contains $S$.
\end{lemma}
\begin{proof}
Since $G$ is connected, for any two elements $u,v$ in $S$ we can choose a finite path
in $G$ with endpoints $u$ and $v$. The subgraph induced by the union of all chosen paths is clearly finite and connected.
\end{proof}

\begin{lemma}\label{lm:connectedremoval}
Any finite connected graph with at least one vertex has a vertex whose removal
makes the remaining graph connected.
\end{lemma}
\begin{proof}
If there is only one vertex, removing it results in the empty graph which is connected.
If there are at least two vertices, let $u$ and $v$ be vertices with maximum
distance. Consider any pair of vertices $x$ and $y$ distinct from $u$. Choose any
shortest path from $x$ to $v$ and any shortest path from $v$ to $y$. Neither of
these paths can be longer than the distance between $u$ and $v$, so they do not
contain $u$. Concatenating these paths shows that $x$ and $y$ are connected in the
graph resulting from removing $u$, so that graph is connected.
\end{proof}

Our main tool in the proof of \cref{th:maindualnormalfield} will be the following
generalization of \cref{th:brownhumphriesnormalfield} to affine transvection groups.
\begin{theo}\label{th:mainnormalfieldaffine}
Suppose $K\ne\bbF_2$, and let $S$ be a spanning subset of $X$
such that $G(S)$ is connected. Then, for any $\alpha\in X^\ast$, the affine transvection group $\transgr{\alpha}{S}$ has at most one non-singleton orbit.
\end{theo}
\begin{proof}
Let $B\subseteq S$ be a basis of $X$. Then, by \cref{lm:infiniteconnectivity}
there is a finite $S'\subseteq S$ containing $B$ such that $G(S')$ is connected;
let us choose $S'$ to be of minimal cardinality with this property.
Then, $\transgr{\alpha}{S'}$ is a subgroup of $\transgr{\alpha}{S}$
and they have the same set of fixed points, namely $\omega^{-1}(-\alpha)$,
so if the theorem holds for $S'$ it holds for $S$ too. Thus, in the following we
may assume that $S$ is finite and that $\Span(S\setminus\{s\})$ is a proper subset
of $X$ for any $s\in S$ such that $G(S\setminus\{s\})$ is connected.
We will use induction on the cardinality of $S$.

Suppose $S$ is a singleton set, $S=\{s\}$. If $\alpha(s)=0$, the group
$\transgr{\alpha}{S}$ acts trivially on $X=Ks$. If $\alpha(s)\ne0$,
$\transgr{\alpha}{S}$ acts transitively on $X$:
Any $x,y\in X$ can be written as
$x=as$ and $y=bs$ for some $a,b\in K$, and putting $k=(b-a)/\alpha(s)$, we obtain
$\trans{\alpha(s)}{s,k}(x)=y$.

In the following we assume that $S$ has more than one element.
Since $G(S)$ is connected and finite, by \cref{lm:connectedremoval} there is an $s\in S$ such that $G(S\setminus\{s\})$ is connected, and, by our earlier assumption,
$Y=\Span(S\setminus\{s\})$ is a proper subspace of $X$.
Let $X_+ := X\setminus X^{\transgr{\alpha}{S}}$.
For each $a\in K$, let $A^a=(as+Y)\cap X_+$,
and partition each $A^a$ as $A^a=A_0^a\cup A_+^a$, where $A_0^a=A^a\cap X^{\transgr{\alpha}{S\setminus\{s\}}}$
and $A_+^a=A^a\setminus X^{\transgr{\alpha}{S\setminus\{s\}}}$. By \cref{lm:groupactionisomorphism},
the permutation group $\restrict{\transgr{\alpha}{S\setminus\{s\}}}{as+Y}$ is isomorphic to the
permutation group $\restrict{\transgr{\alpha+\omega(as)}{S\setminus\{s\}}}{Y}$, which, by induction,
has at most one non-singleton orbit, so $\transgr{\alpha}{S\setminus\{s\}}$, acts transitively on
each $A_+^a$. For any $as+y\in A_0^a$, we have $p:=\omega(as+y,s)+\alpha(s)\ne0$, so
for any $b\ne a$ in $K$, we have $\trans{\alpha(s)}{s,(b-a)/p}(x) = bs + y$. Since
$G(S)$ is connected and $S$ has at least two elements, there is a
$t\in S\setminus\{s\}$ such that $\omega(s,t)\ne0$. Since
$as+y\in A_0^a\subseteq X^{\transgr{\alpha}{S\setminus\{s\}}}$,
we have $\omega(as+y,t)+\alpha(t)=0$, and it follows that $\omega(bs+y,t)+\alpha(t)\ne0$.
Hence, $bs+y$ belongs to $A_+^b$.

We have shown that every element in $A_0^a$ belongs to the same $\transgr{\alpha}{S}$-orbit as the elements in $A_+^b$ with $b\ne a$, and since
$K$ has more than two elements, it follows that
all of $X_+$ belongs to the same orbit.
\end{proof}

Our first main result is the following dual analogue to \cref{th:brownhumphriesFtwo}.
\maindualnormalfield
\begin{proof}
Let $\theta_\alpha$ be the affine mapping from $X$ to $X^\ast$ defined by
$\theta_\alpha(x) = \omega(x) + \alpha$.
Note that $\theta_\alpha^{-1}(\alpha)=\ker\omega$, which is nonempty.
By \cref{lm:pressingstolamps}, $\alpha$ and $\beta$ belong to the same
$\transgrsimple{S}^\ast$-orbit if and only if $\theta_\alpha^{-1}(\beta)$ is nonempty too, and belongs to the same $\restrict{\transgr{\alpha}{S}}{X/\ker\omega}$-orbit as $\theta_\alpha^{-1}(\alpha)$. Clearly,
$\theta_\alpha^{-1}(\beta)$ is nonempty
if and only if $\beta-\alpha$ belongs to $\im\omega$.
Note that, for any $x\in\theta_\alpha^{-1}(\alpha)=\ker\omega$ and $y\in\theta_\alpha^{-1}(\beta)$, both $\omega(x)+\alpha=\alpha$ and $\omega(y)+\alpha=\beta$
are nonzero, so neither $x$ nor $y$ is a fixed point of $\transgr{\alpha}{S}$.
Hence, by \cref{th:mainnormalfieldaffine}, $x$ and $y$
belong to the same $\transgr{\alpha}{S}$-orbit.
We conclude that $\alpha$ and $\beta$ belong to the same
$\transgrsimple{S}^\ast$-orbit if and only if $\beta-\alpha\in\im\omega$.
\end{proof}

\section{\texorpdfstring{The case $K=\bbF_2$ and $S$ is a basis of orthogonal type}{Two-element field and basis of orthogonal type}}\label{sec:orthogonal}
The goal of this section is to prove our second main result, \cref{th:maindualmodtwoorthogonal}. The idea is the same as for the case where $K\ne\bbF_2$, namely to first consider orbits of affine transvection groups.
Those orbits are described by \cref{th:mainmodulotwoaffine} below, the important
part of which was proved already in 2000 by Shapiro
et~al.~\cite[Th.~7.2]{ShapiroEtAl2000}.
Their proof relies on \cref{th:brownhumphriesFtwo} (in the form
of their Lemma~7.7 which is Lemma~3.4 in \cite{ShapiroEtAl1998}, which in turn
depends on Theorem~3.5 in \cite{Janssen1983}).
Our proof is self-contained and perhaps conceptually simpler, so we
hope it has some independent value.

In this section we assume that $X$ is a finite-dimensional vector space over
$K=\bbF_2$ equipped with an alternating bilinear form $\omega$, and that
$S$ is a basis of $X$. Note that $X$, $\omega$ and $S$ are recoverable from the graph $G(S)$.

Though we will not use this terminology here, it might be useful to adopt the
intuition from the lit-only sigma game and think about
an element $x\in X$ as a pressing configuration on the vertices $S$ of the graph
$G(S)$, such that a vertex $s\in S$ is pressed if the $s$-coordinate of $x$ is
one. The value $Q_S(x)$ of the quadratic form equals the number of pressed vertices
plus the number of edges between them modulo two, that is, the Euler
characteristic of the subgraph induced by pressed vertices.
An element $\beta\in X^\ast$ can be thought of as a lamp configuration on the
vertices $S$ of $G(S)$, such that a vertex $s\in S$ is lit if and only if
$\beta(s)=1$.
In this setting, \cref{lm:pressingstolamps} can be interpreted as follows.
Let each pressing configuration $x$ automatically yield the lamp configuration $\omega(x)+\alpha$. Applying an affine transvection
$\trans{\alpha(s)}{s}$ to $x$
has no effect if $s$ is not lit, and if $s$ is lit it has
the effect of toggling the button at $s$ and toggling the lamp at each neighbor of
$s$. For the lamp configuration this is equivalent to applying the dual tranvection $\transsimple{s}^\ast$.

It is known that the transvection group
$\transgrsimple{S}$ preserves the quadratic form $Q_S$; see
e.g.~\cite{BrownHumphriesII1986}.
The following proposition generalizes this result to affine transvection groups.
\begin{prop}\label{pr:Qplusalphainvariant}
Let $S$ be a basis of $X$. Then, $\transgr{\alpha}{S}$ preserves $Q_S+\alpha$.
\end{prop}
\begin{proof}
Take any $s\in S$ and any $x\in X$.
Since $Q_S(x+s)=Q_S(x)+Q_S(s)+\omega(x,s)$, we have
\begin{align*}
\trans{\alpha(s)}{s}(x) &= x+\bigl(\omega(x,s)+\alpha(s)\bigr)s\\
&=x+\bigl(Q_S(x+s)+Q_S(x)+Q_S(s)+\alpha(s)\bigr)s\\
&=x+\bigl(Q_S(s)+(Q_S+\alpha)(x+s)+(Q_S+\alpha)(x)\bigr)s\\
&=x+\bigl(1+(Q_S+\alpha)(x+s)+(Q_S+\alpha)(x)\bigr)s.
\end{align*}
If $(Q_S+\alpha)(x+s)+(Q_S+\alpha)(x)=1$, we have $\trans{\alpha(s)}{s}(x)=x$,
and $(Q_S+\alpha)(\trans{\alpha(s)}{s}(x))=(Q_S+\alpha)(x)$ holds trivially.
If $(Q_S+\alpha)(x+s)+(Q_S+\alpha)(x)=0$, we have
$\trans{\alpha(s)}{s}(x)=x+s$, so $(Q_S+\alpha)(\trans{\alpha(s)}{s}(x))=(Q_S+\alpha)(x+s)
=(Q_S+\alpha)(x)$.
\end{proof}

Next, we want to show that t-equivalent bases behave the same with regard to
quadratic forms and orbits. To this end, we need the following lemma.
\begin{lemma}\label{lm:tripptrapptrull}
For any $s,t\in X$ with $\omega(s,t)=1$ and any $a,b\in K$,
it holds that $\trans{a}{s}\circ \trans{b}{t}\circ \trans{a}{s}=\trans{a+b}{s+t}$.
\end{lemma}
\begin{proof}
This is a tedious but straightforward matter of applying the definitions
and using that $\omega(s,s)=0$ and $\omega(s,t)=1$:
\begin{align*}
\trans{a}{s}(x) &= x + [a+\omega(s,x)]s,\\
\trans{b}{t}(\trans{a}{s}(x)) &= x + [a+\omega(s,x)]s + [b+\omega(t,x + [a+\omega(s,x)]s)]t\\
&=x + [a+\omega(s,x)]s + [a+b+\omega(s+t,x)]t,\\
\trans{a}{s}(\trans{b}{t}(\trans{a}{s}(x)))
&= x + [a+\omega(s,x)]s + [a+b+\omega(s+t,x)]t \\
&\quad + [a+\omega(s,x+[a+\omega(s,x)]s + [a+b+\omega(s+t,x)]t)]s\\
&=x + [a+b+\omega(s+t,x)](s+t) \\
&= \trans{a+b}{s+t}(x).
\end{align*}
\end{proof}

\begin{prop}\label{pr:equiv}
If two bases $S$ and $S'$ of $X$ are t-equivalent, then $Q_S=Q_{S'}$ and
$\transgr{\alpha}{S}=\transgr{\alpha}{S'}$ for any $\alpha\in X^\ast$.
\end{prop}
\begin{proof}
It is enough to check the case where $S'$ can be obtained from $S$ by a single
t-equivalence step, that is, by replacing $t$ with $s+t$ for some $s,t\in S$ with
$\omega(s,t)=1$.

By \cref{lm:tripptrapptrull},
$\trans{\alpha(s+t)}{s+t}=\trans{\alpha(s)}{s}\circ \trans{\alpha(t)}{t}\circ \trans{\alpha(s)}{s}\in\transgr{\alpha}{S}$, so $\transgr{\alpha}{S'}\subseteq\transgr{\alpha}{S}$.
Also, $\trans{\alpha(t)}{t}=\trans{\alpha(s)}{s}\circ\trans{\alpha(s+t)}{s+t}\circ\trans{\alpha(s)}{s}\in\transgr{\alpha}{S'}$, so $\transgr{\alpha}{S}\subseteq\transgr{\alpha}{S'}$
and we conclude that $\transgr{\alpha}{S}=\transgr{\alpha}{S'}$.
To see that $Q_S=Q_{S'}$, it is enough to check that $Q_S(s')=1$ for any
$s'\in S'$. But the only element in $S'$ that does not belong to $S$ is $s+t$, and
$Q_S(s+t)=Q_S(s)+Q_S(t)+\omega(s,t)=1+1+1=1$.
\end{proof}

\begin{defi}
We say that a basis $S$ of $X$, and the corresponding graph $G(S)$, is \emph{connecting}
if, for any $\alpha\in X^\ast$, any two nonfixed points of
$\transgr{\alpha}{S}$ belong to the same orbit if and only if
they have the same $(Q_S+\alpha)$-value.
\end{defi}
Note that, by \cref{pr:Qplusalphainvariant}, the ``only if'' part of the definition
holds for any basis $S$.

\begin{prop}\label{pr:Esixconnecting}
The graph $E_6=%
\begin{tikzpicture}[scale=0.5]
\draw (0,0)--(4,0);
\draw (2,0)--(2,1);
\filldraw (0,0) circle (3pt);
\filldraw (1,0) circle (3pt);
\filldraw (2,0) circle (3pt);
\filldraw (3,0) circle (3pt);
\filldraw (4,0) circle (3pt);
\filldraw (2,1) circle (3pt);
\end{tikzpicture}$
is connecting.
\end{prop}
\begin{proof}
Let $S$ be the vertex set of $E_6$, and let
$X=\generate{S}$ with $\omega(s,t)=1$ for $s,t\in S$ if and only if $s$ and $t$ are
neighbors in the graph.
It is easy to verify by hand (and it also follows from \cref{th:brownhumphriesFtwo})
that, under the action of $\transgrsimple{S}$
on $X$, two nonfixed elements belong to the same orbit if and only
if they have the same $Q_S$-value. It is also easy to check that
the bilinear form $\omega$ is nondegenerate.
Let $\alpha$ be any element in $X^\ast$ and let
$x$ be the unique element in $X$ such that $\omega(x) = \alpha$.

Let $\psi:X\rightarrow X$ be the mapping that adds $x$.
By \cref{lm:groupactionisomorphism}, $\psi$ induces a group isomorphism between
$\transgrsimple{S}$ and
$\restrict{\transgr{\alpha}{S}}{x+X}=\transgr{\alpha}{S}$, and
for any $y\in X$ we have $(Q_S+\alpha)(\psi(y)) = (Q_S+\alpha)(y+x) = Q_S(y) + Q_S(x)$, so
it follows that, under the action of $\transgr{\alpha}{S}$ on $X$, two nonfixed elements belong to the same orbit if and only if they have the same $(Q_S+\alpha)$-value.
\end{proof}

We want to show that being connecting is a monotone graph property, so that
if a connected graph contains a connecting graph as an induced subgraph, then
the larger graph would be connecting too. To make the induction step work,
however, we need an additional property:

\begin{defi}
We say that a basis $S$ of $X$, and the corresponding graph $G(S)$,
is \emph{nice} if, for any
$\chi\in \bbF_2$ and any $\alpha,\beta\in X^\ast$, there is an
$x\in X$ such that $(Q_S+\alpha)(x)=(Q_S+\beta)(x)=\chi$ and
$\omega(x)\not\in\{\alpha,\beta\}$.
\end{defi}

\begin{prop}\label{pr:Afournice}
The graph \twoDisjointEdges\ is nice.
\end{prop}
\begin{proof}
This can of course be verified easily with a computer, but we prefer to give a
human proof.

Let $X$ be the two-dimensional vector space with basis $S=\{s_1,s_2\}$ and
symplectic form $\omega_X(s_1,s_2)=1$, and consider any $\alpha,\beta\in X^\ast$.
We will express $\alpha$ and $\beta$ together by the matrix
$A=\amat{\alpha(s_1)}{\alpha(s_2)}{\beta(s_1)}{\beta(s_2)}$.
Now, construct a set
$M(A)$ as follows. For any $x\in X$, let $\bfm(x)$ be the vector
$\mvec{(Q_S+\alpha)(x)}{(Q_S+\beta)(x)}$
in $\bbF_2^2$, and mark the upper and lower entry by a star if
$\omega_X(x)\ne\alpha$ and $\omega_X(x)\ne\beta$, respectively.
Let $M(A)=\{\bfm(x):x\in X\}$.

For instance, if $A=\amat0001$ then
$\bfm(0)=\mvec{0}{\demark{0}}$, 
$\bfm(s_1)=\mvec{\demark{1}}{1}$ and
$\bfm(s_2)=\bfm(s_1+s_2)=\mvec{\demark{1}}{\demark{0}}$,
so $M(A)=\{
\mvec{0}{\demark{0}},
\mvec{\demark{1}}{1},
\mvec{\demark{1}}{\demark{0}}
\}$.
In \cref{tb:M},
\begin{table}
\[
\setlength{\cellspacetoplimit}{5pt}
\setlength{\cellspacebottomlimit}{5pt}
\begin{tabular}{>{$}Sr<{$}|>{$}Sl<{$}}
A & M(A) \\ \hline
\Amat0000
& \left\{
\Mvec{0}{0}, \Mvec{\demark{1}}{\demark{1}}
\right\}
\\ \hline
\Amat0010, \Amat0001, \Amat0011
& \left\{
\Mvec{0}{\demark{0}}, \Mvec{\demark{1}}{1}, \Mvec{\demark{1}}{\demark{0}}
\right\}
\\ \hline
\Amat1000, \Amat0100, \Amat1100
& \left\{
\Mvec{\demark{0}}{0}, \Mvec{1}{\demark{1}}, \Mvec{\demark{0}}{\demark{1}}
\right\}
\\ \hline
\Amat1010, \Amat0101, \Amat1111
& \left\{
\Mvec{\demark{0}}{\demark{0}}, \Mvec{1}{1}
\right\}
\\ \hline
\Amat1001, \Amat0110, \Amat1110, \Amat1101, \Amat0111, \Amat1011
& \left\{
\Mvec{\demark{0}}{\demark{0}}, \Mvec{1}{\demark{0}}, \Mvec{\demark{0}}{1}
\right\}
\end{tabular}
\]
\caption[$M(A)$ for all possible $A$.]{$M(A)$ for all possible $A$. Note that $M(A)$ is invariant
under reordering of the columns of $A$, and reordering the rows of $A$ just
reorders the rows of $\bfm(x)$ correspondingly, so when generating the table we
only had to consider the $A$-values
$\amat0000$, $\amat0010$, $\amat0011$, $\amat1010$,
$\amat1001$, $\amat1110$ and $\amat1111$.}\label{tb:M}
\end{table}
we have computed $M(A)$ for all possible $A$.

The vector space corresponding to the graph \twoDisjointEdges\ is
the direct sum of two
copies of $X$. To show that this graph is nice we need to check that,
for any $\chi\in\bbF_2$ and for any pair of $A$-values $A_1$ and $A_2$ in the table, there are an
$\bfm_1$ in $M(A_1)$ and an $\bfm_2$ in $M(A_2)$ such that
$\bfm_1+\bfm_2=\mvec{\demark{\chi}}{\demark{\chi}}$,
where an entry in the sum is defined to be marked by a star if at least one of the corresponding entries in the summands has a star.
In \cref{fig:mvecrelations}, we have drawn solid edges between some pair
of $\bfm$-vectors whose sum is $\mvec{\demark{1}}{\demark{1}}$
and dashed edges between some pair of $\bfm$-vectors whose sum is $\mvec{\demark{0}}{\demark{0}}$. We see that
any $M(A_1)$ and $M(A_2)$ are connected by both a solid and a dashed edge.
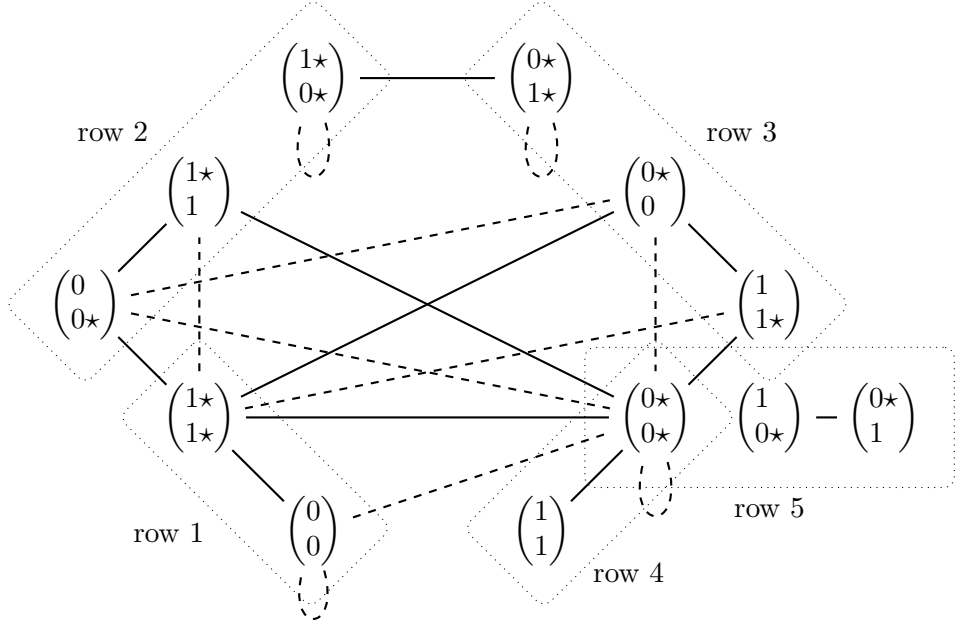
\begin{figure}
\begin{tikzpicture}[scale=0.15, every node/.style={}, every fit/.style={rounded corners,draw,dotted,inner sep=9pt}]
\begin{scope}[every node/.style={circle, inner sep={-1pt}, outer sep=0}]
\node (A) at (-10,-20) {$\Mvec{0}{0}$};
\node (B) at (-20,-10) {$\Mvec{\demark{1}}{\demark{1}}$};
\node (C) at (-30,0) {$\Mvec{0}{\demark{0}}$};
\node (D) at (-20,10) {$\Mvec{\demark{1}}{1}$};
\node (E) at (-10,20) {$\Mvec{\demark{1}}{\demark{0}}$};
\node (F) at (10,20) {$\Mvec{\demark{0}}{\demark{1}}$};
\node (G) at (20,10) {$\Mvec{\demark{0}}{0}$};
\node (H) at (30,0) {$\Mvec{1}{\demark{1}}$};
\node (I) at (20,-10) {$\Mvec{\demark{0}}{\demark{0}}$};
\node (J) at (10,-20) {$\Mvec{1}{1}$};
\node (K) at (30,-10) {$\Mvec{1}{\demark{0}}$};
\node (L) at (40,-10) {$\Mvec{\demark{0}}{1}$};
\end{scope}
\begin{scope}[thick]
\path[-] (A) edge (B);
\path[-] (B) edge (C);
\path[-] (C) edge (D);
\path[-] (E) edge (F);
\path[-] (B) edge (G);
\path[-] (G) edge (H);
\path[-] (H) edge (I);
\path[-] (I) edge (J);
\path[-] (D) edge (I);
\path[-] (B) edge (I);
\path[-] (K) edge (L);
\end{scope}
\begin{scope}[dashed, thick, every loop/.style={}]
     \draw (A) -- (I);
     \draw (B) -- (H);
     \draw (B) -- (D);
     \draw (C) -- (G);
     \draw (C) -- (I);
     \draw (G) -- (I);
     \draw (A) edge [loop below] (A);
     \draw (E) edge [loop below] (E);
     \draw (F) edge [loop below] (F);
     \draw (I) edge [loop below] (I);
\end{scope}
\node [rotate fit=45, fit=(A)(B), label=left:row 1] {};
\node [rotate fit=45, fit=(C)(D)(E), label=above:row 2] {};
\node [rotate fit=45, fit=(F)(G)(H), label=right:row 3] {};
\node [rotate fit=45, fit=(I)(J), label=below left:\ \ \ \ \ \ row 4] {};
\node [rotate fit=0, fit=(I)(K)(L), label=below:row 5] {};
\end{tikzpicture}
\caption[Some relations of the $\bfm$-vectors.]{Some relations of the $\bfm$-vectors. Two elements connected
by a solid edge sum to $\mvec{\demark{1}}{\demark{1}}$, and two
elements connected by a dashed edge sum to
$\mvec{\demark{0}}{\demark{0}}$. The framed sets correspond to the
rows of \cref{tb:M}.}\label{fig:mvecrelations}
\end{figure}
\end{proof}

Next, we show that niceness is a monotone graph property.
\begin{prop}\label{pr:nicenessmonotone}
Let $S$ be a basis of $X$ and let $T$ be a subset of $S$.
Suppose that $T$ is a nice basis of its span.
Then, $S$ is nice too.
\end{prop}
\begin{proof}
Let $Y=\Span{T}$ and take any $\alpha,\beta\in X^\ast$ and any $\chi\in\bbF_2$.
Since $T$ is nice, there is a $y\in Y$ such that
$(\restrict{Q_S}{Y}+\restrict{\alpha}{Y})(y)=(\restrict{Q_S}{Y}+\restrict{\beta}{Y})(y)=\chi$ and
$\restrict{\omega_X(y)}{Y}\not\in\{\restrict{\alpha}{Y},\restrict{\beta}{Y}\}$.
It follows that $(Q_S+\alpha)(y)=(Q_S+\beta)(y)=\chi$ and
$\omega_X(y)\not\in\{\alpha,\beta\}$.
Hence, $S$ is a nice basis.
\end{proof}

Before we are ready for the induction proof, we need one more proposition.
\begin{prop}\label{pr:niceisnice}
Suppose $S$ is a nice basis of $X$. Let $\chi,\psi\in\bbF_2$ and let $\alpha$ and $\beta$ be distinct elements of $X^\ast$.
Then there is an $x\in X$ such that $(Q_S+\alpha)(x)=\psi$ and $(Q_S+\beta)(x)=\chi$ and $\omega(x)\not\in\{\alpha,\beta\}$.
\end{prop}
\begin{proof}
If $\chi=\psi$, the result follows directly from the definition of niceness, so let us
assume that $\psi=\chi+1$

Take any $y\in X$ such that $(\alpha+\beta)(y)=1$. Put $\chi'=\chi+(Q_S+\beta)(y)$, $\alpha'=\alpha+\omega(y)$ and $\beta'=\beta+\omega(y)$. Since $S$ is nice there is an
$x'\in X$ such that $(Q_S+\alpha')(x')=(Q_S+\beta')(x')=\chi'$ and
$\omega(x')\not\in\{\alpha',\beta'\}$. Set $x:=x'+y$. This means that
\begin{align*}
(Q_S+\alpha)(x) + (Q_S+\alpha)(y) &= Q_S(x'+y) + Q_S(y) + \alpha(x') \\
&= Q_S(x') + \omega(x',y) + \alpha(x') \\
&=(Q_S+\alpha')(x') = \chi',
\end{align*}
and, by a symmetric argument, $(Q_S+\beta)(x) + (Q_S+\beta)(y) = \chi'$ as well.
Hence, $(Q_S+\alpha)(x) = \chi' + (Q_S+\alpha)(y) = \chi + (\alpha+\beta)(y)
= \chi + 1 = \psi$ and 
$(Q_S+\beta)(x) = \chi' + (Q_S+\beta)(y) = \chi$.
Also,
$\omega(x')+\omega(y)+\alpha=\omega(x)+\alpha$ and 
$\omega(x')+\omega(y)+\beta=\omega(x)+\beta$,
so $\omega(x)\not\in\{\alpha,\beta\}$.
\end{proof}

Finally, we have everything we need to show by induction
that being connecting and nice
is a monotone graph property for connected graphs.
\begin{prop}\label{pr:induction}
Let $S$ be a basis of $X$ such that $G(S)$ is connected,
and let $s\in S$. Suppose
$S$ has at least two elements and $S\setminus\{s\}$ is a nice and connecting basis of its span.
Then, $S$ is nice and connecting.
\end{prop}
\begin{proof}
Take any $\alpha\in X^\ast$ and any $\chi\in\bbF_2$.
Let
\[
X_+ = \{x\in X\setminus X^{\transgr{\alpha}{S}}\,:\,(Q_S+\alpha)(x)=\chi\}
\]
be the set of nonfixed points with $(Q_S+\alpha)$-value $\chi$.
We must show that $\transgr{\alpha}{S}$ acts transitively on $X_+$.

Since $G(S)$ is connected there is a $t$ in $S\setminus\{s\}$ such that
$\omega(s,t)=1$.

Take any $x\in X_+$. If $x$ is a fixed point of $\transgr{\alpha}{S\setminus\{s\}}$,
we have $\trans{\alpha(s)}{s}(x)=x+s$ which is not a fixed point of
$\transgr{\alpha}{S\setminus\{s\}}$ since
$(\omega(x+s)+\alpha)(t)=\omega(s,t)+(\omega(x)+\alpha)(t)=1+0=1$.
Hence it remains only to
show that $\transgr{\alpha}{S}$ acts transitively on
$X_+\setminus X^{\transgr{\alpha}{S\setminus\{s\}}}$. This set can be divided
into two parts: $\calA$, consisting of those $x$ that belong to $\Span(S\setminus\{s\})$,
and $\calB$, consisting of those $x$ that belong to $s+\Span(S\setminus\{s\})$.
Since $S\setminus\{s\}$ is a nice and connecting basis of its span,
$\transgr{\alpha}{S\setminus\{s\}}$ acts transitively on
$\calA$. By \cref{lm:groupactionisomorphism} (with $x$ set to $s$), it acts transitively on $\calB$ too.

Thus, it suffices to show that there is some element in $\calA$ that belong to the same $\transgr{\alpha}{S}$-orbit as some element in $\calB$.

$\omega(s)$ is not identically zero on $S\setminus\{s\}$ since
$\omega(s,t)=1$, so we can apply \cref{pr:niceisnice} on $\chi$ and $\psi=\chi+(Q_S+\alpha)(s)$ and the linear forms $\alpha$ and $\alpha+\omega(s)$. Hence, there is an $x\in\Span(S\setminus\{s\})$ such that
$(Q_S+\omega(s)+\alpha)(x)=\chi+(Q_S+\alpha)(s)$ and $\chi=(Q_S+\alpha)(x)$
and neither $\omega(x)+\alpha$ nor $\omega(x)+\omega(s)+\alpha$ is identically zero on
$S\setminus\{s\}$. It follows that $x$ belongs to $\calA$. We have
$(\omega(x)+\alpha)(s)=(Q_S+\omega(s)+\alpha)(x) + (Q_S+\alpha)(x) + \alpha(s)=
\chi + (Q_S+\alpha)(s) + \chi + \alpha(s) = Q_S(s)=1$, so
$\trans{\alpha(s)}{s}(x)=x+s$ which belongs to $\calB$ since $\omega(x+s)+\alpha$
is not identically zero on $S\setminus\{s\}$.
\end{proof}

The following theorem generalizes \cref{th:brownhumphriesFtwo} to
affine transvection groups. As noted above, part of it was proved by Shapiro et.~al~\cite[Th.~7.2]{ShapiroEtAl2000}.
\begin{theo}\label{th:mainmodulotwoaffine}
Any basis of orthogonal type is nice and connecting.
\end{theo}
\begin{proof}
By \cref{pr:equiv}, we may assume that
$G(S)$ contains $E_6$ as an induced subgraph.
By \cref{pr:Esixconnecting,pr:Afournice,pr:nicenessmonotone}, $E_6$ is both nice and connecting,
so by \cref{pr:induction} it follows by induction that $S$ is nice and connecting.
\end{proof}

As a final ingredient in the proof of \cref{th:maindualmodtwoorthogonal} we need
the following lemma, which shows that $Q_S+\alpha$ is either constant on each coset in $X/\ker\omega$ or nonconstant on each coset.
\begin{lemma}\label{lm:qorbits}
If there is an $x_0\in\ker\omega$ such that
$Q_S(x_0)\ne\alpha(x_0)$, then $(Q_S+\alpha)(p)=\{0,1\}$
for any coset $p\in X/\ker\omega$.
Otherwise, $Q_S+\alpha$ is constant on each coset.
\end{lemma}
\begin{proof}
Take any $p\in X/\ker\omega$ and write $p=x+\ker\omega$ for some $x\in X$.
If there is an $x_0\in\ker\omega$ such that $Q_S(x_0)\ne\alpha(x_0)$,
then $x+x_0$ belongs to $p$ and
\[
(Q_S+\alpha)(x+x_0)=(Q_S+\alpha)(x) + (Q_S+\alpha)(x_0) + \omega(x,x_0)
=(Q_S+\alpha)(x) + 1 + 0.
\]
We conclude that $(Q_S+\alpha)(p)=\{0,1\}$.

If $Q_S(x_0)=\alpha(x_0)=0$ for any $x_0\in\ker\omega$, then
for any $x_0\in\ker\omega$,
\[
(Q_S+\alpha)(x+x_0)=(Q_S+\alpha)(x) + (Q_S+\alpha)(x_0) + \omega(x,x_0)
=(Q_S+\alpha)(x) + 0 + 0,
\]
so $Q_S+\alpha$ is constant on $p$.
\end{proof}

At last, we are ready to prove our second main result, which is a dual analogue
of \cref{th:brownhumphriesFtwo}.
\maindualmodtwoorthogonal
\begin{proof}
The first part of the proof is identical to first part of the proof of
\cref{th:maindualnormalfield}.

Let $\theta_\alpha$ be the affine mapping from $X$ to $X^\ast$ defined by
$\theta_\alpha(x) = \omega(x) + \alpha$.
Note that $\theta_\alpha^{-1}(\alpha)=\ker\omega$, which is nonempty.
By \cref{lm:pressingstolamps}, $\alpha$ and $\beta$ belong to the same
$\transgrsimple{S}^\ast$-orbit if and only if $\theta_\alpha^{-1}(\beta)$ is nonempty too, and belongs to the same $\restrict{\transgr{\alpha}{S}}{X/\ker\omega}$-orbit as $\theta_\alpha^{-1}(\alpha)$. Clearly,
$\theta_\alpha^{-1}(\beta)$ is nonempty
if and only if $\beta-\alpha$ belongs to $\im\omega$.
Note that, for any $x\in\theta_\alpha^{-1}(\alpha)=\ker\omega$ and $y\in\theta_\alpha^{-1}(\beta)$, both $\omega(x)+\alpha=\alpha$ and $\omega(y)+\alpha=\beta$
are nonzero, so neither $x$ nor $y$ is a fixed point of $\transgr{\alpha}{S}$.
By \cref{th:mainmodulotwoaffine}, $S$ is a connecting basis.
Hence, $\theta_\alpha^{-1}(\alpha)$ and $\theta_\alpha^{-1}(\beta)$ belong to the same orbit if and only if
$(Q_S+\alpha)(\theta_\alpha^{-1}(\alpha))\cap(Q_S+\alpha)(\theta_\alpha^{-1}(\beta))\ne\emptyset$. Since $0\in\theta_\alpha^{-1}(\alpha)$ and $(Q_S+\alpha)(0)=0$, by \cref{lm:qorbits} this happens if and only if there
is an $x\in\theta_\alpha^{-1}(\beta)$ such that $(Q_S+\alpha)(x)=0$.
\end{proof}

\section{\texorpdfstring{The case $K=\bbF_2$ and $S$ is a basis not of orthogonal type}{Bases not of orthogonal type}}\label{sec:linegraphs}
Suppose $K=\bbF_2$, and let $S$ be a basis of $X$ not of orthogonal type such that $G(S)$ is connected. By \cref{th:cuypers}, the graph $G(S)$ is the
line graph of some connected multigraph $G=(V,E)$, so we can identify $S$ with $E$,
and we adopt the notation from \cref{sec:results}.

Our approach will be very similar to the one taken by Wu~\cite{Wu2009}.
The difference is that Wu considers only line graphs of \emph{simple} graphs and focuses on the size of orbits rather than trying to answer \cref{qu:main}.

As in \cref{sec:orthogonal}, it might help the intuition to interpret the
situation in terms of buttons and lamps.
We can think of an element $y\in\generate{V}$ as a pressing configuration on
the vertices $V$, such that a vertex $v\in V$ is pressed if the $v$-coordinate
of $y$ is one. 
An element $\beta\in X^\ast$ can be thought of as a lamp configuration on the
edges $E$ of $G$, such that an edge $e\in E$ is lit if and only if
$\beta(e)=1$.

Note that, since $K=\bbF_2$, the inner product $\omega_{\generate{V}}$ is
a skew-symmetric bilinear form.
\begin{prop}\label{pr:deltapartialomega}
$\partial$ preserves the bilinear form, and $\delta\circ\partial=\omega$.
\end{prop}
\begin{proof}
For any edges $e,e'\in E$, it holds that
$\omega_{\generate{V}}(\partial(e),\partial(e'))$ is 1 if and only if
$e$ and $e'$ have exactly one common endpoint, which happens if and only if
$\omega(e,e')=1$. By bilinearity, it follows that
$\omega_{\generate{V}}(\partial(x),\partial(y))=\omega(x,y)$ for any $x,y\in X$.

Now, for any $x,y\in X$ we have
$(\delta\circ\partial)(x)(y)=\omega_{\generate{V}}(\partial(x),\partial(y))
=\omega(x,y)$.
\end{proof}

Let $E_\mathrm{span}\subseteq E$ be the edge set of a spanning tree in $G$.
Also, let $\ones:=\sum_{v\in V}v$ denote the configuration with all vertices pressed.
\begin{lemma}\label{lm:directsum}
$X^\ast=\im\delta\oplus E_\mathrm{span}^0$,
where $E_\mathrm{span}^0$ is the annihilator of $E_\mathrm{span}$.
\end{lemma}
\begin{proof}
Take any $\beta\in X^\ast$.
Pick any vertex in $V$ as the \emph{root} of the tree $E_\mathrm{span}$,
and define $y\in\generate{V}$ by, for each vertex $v\in V$,
letting the $v$-coordinate of $y$ be the
sum of the $\beta$-values of all edges along the unique path from the
root to $v$. Then $\delta(y)=\omega_{\generate{V}}(y)\circ\partial$ coincides with $\beta$ on the set $E_\mathrm{span}$.
The difference $\beta - \delta(y)$ belongs to $E_\mathrm{span}^0$, and we conclude
that $X^\ast=\im\delta + E_\mathrm{span}^0$.

To show that $\im\delta\cap E_\mathrm{span}^0=\{0\}$,
suppose $\beta\in\im\delta\cap E_\mathrm{span}^0$. Then, $\beta=\delta(y)
=\omega_{\generate{V}}(y)\circ\partial$ for some $y\in\generate{V}$, and
$\omega_{\generate{V}}(y,\partial(e))=(\omega_{\generate{V}}(y)\circ\partial)(e)=\beta(e)=0$ for any $e\in E_\mathrm{span}$. Hence, $y$ 
has the same coordinates at the endpoints of
each edge of the spanning tree $E_\mathrm{span}$,
and we conclude that $y$ is either $0$ or $\ones$ and that $\beta=\omega_{\generate{V}}(y)\circ\partial=0$.
\end{proof}

Recall the definitions from \cref{sec:notation}.
\begin{lemma}\label{lm:pressingstolampslinegraph}
Let $\delta_\alpha$ be the mapping from $\generate{V}$ to $X^\ast$ defined by
$\delta_\alpha(y) = \delta(y) + \alpha$.
Then
$\delta_\alpha$ induces a group isomorphism from
$\restrict{\transgr{\alpha}{E,\partial}}{\generate{V}/\ker\delta}$ to
$\restrict{\transgrsimple{E}^\ast}{\im\delta_\alpha}$
\end{lemma}
\begin{proof}
By \cref{pr:deltapartialomega}, $\partial$ preserves the bilinear form,
so the lemma follows from \cref{lm:pressingstolampsmaster} with $Y=\generate{V}$
and $\phi=\partial$.
\end{proof}
The lemma can be interpreted in terms of pressings and lamps:
Let each pressing configuration $y\in\generate{V}$ automatically yield the lamp configuration $\delta(y)+\alpha$. Applying an affine transvection
$\trans{\alpha(s)}{\partial(e)}$ to $y$
has no effect if $e$ is not lit, and if $e$ is lit it has
the effect of toggling the buttons at the endpoints of $e$ and toggling the
lamp at each edge that has exactly one endpoint in common with $e$. For the lamp configuration this is equivalent to applying the dual tranvection $\transsimple{e}^\ast$.

For an element $y$ in $\generate{V}$, recall that $d_0(y)$ and $d_1(y)$ denote the number of
zero and one coordinates of $y$ in the basis $V$, respectively.
and that $d(y)=\min\{d_0(y),d_1(y)\}$.
\begin{prop}\label{pr:treeorbits}
If $\alpha=0$, then two elements $y,z\in\generate{V}$ belong to the same orbit
of $\transgr{\alpha}{E,\partial}$
if and only if $d_1(x)=d_1(y)$. If $\alpha\ne0$ and $\alpha\in E_\mathrm{span}^0$,
then two elements
$y,z\in\generate{V}$ belong to the same orbit
if and only if $d_1(x)$ and $d_1(y)$ have the same parity.
\end{prop}
\begin{proof}
Let us refer to applying $\trans{\alpha(e)}{\partial(e)}$ as ``playing'' the edge
$e\in E$.

If $\alpha=0$, the effect of playing an edge is simply to swap
the coordinates of its endpoints. Since $G$ is connected,
we can redistribute the coordinate values arbitrarily among the vertices by
playing, so $y,z\in\generate{V}$ belong to the same
orbit if and only if $d_1(x)=d_1(y)$.

Now suppose $\alpha\ne0$ and $\alpha\in E_\mathrm{span}^0$. Playing an edge
in the spanning tree $E_\mathrm{span}$ still has the effect of swapping the
coordinates of its endpoints, so we may still redistribute the coordinate values
arbitrarily among the vertices by playing. However, there is at least
one edge $e$ outside $E_\mathrm{span}$ such that $\alpha(e)=1$, and by playing that
edge we can increase or decrease the number of one-coordinates by two:
To increase by two we first redistribute the coordinate values so that the
coordinates are zero at both endpoints of $e$ before playing it.
To decrease by two we first redistribute the coordinates so that the coordinates are
one at both endpoints of $e$ before playing it.
\end{proof}

Note that $\ker\delta=\{0,\ones\}$.
Since $d(y)=d(y+\ones)$, the function $d$ is well defined
on $\generate{V}/\ker\delta$ as well.
\begin{lemma}\label{lm:sameorbitfromd}
If $\alpha=0$, then two elements $p,q\in\generate{V}/\ker\delta$ belong to the same orbit of $\restrict{\transgr{\alpha}{E,\partial}}{\generate{V}/\ker\delta}$ if and only if $d(p)=d(q)$.
If $\alpha\ne0$ and $\alpha\in E_\mathrm{span}^0$, then two elements
$p,q\in\generate{V}/\ker\delta$ belong to the same orbit
if and only if $d(p)-d(q)$ is even or $\#V$ is odd.
\end{lemma}
\begin{proof}
Let $y\in p$ and $z\in q$. Then $p$ and $q$ belong to the same orbit if and only if
(a) $y$ and $z$ belong to the same orbit or (b) $y$ and $z+\ones$ belong to the same orbit.

First suppose $\alpha=0$. Then, by \cref{pr:treeorbits}, case (a) happens if and only if
$d_1(y)=d_1(z)$ and case (b) happens if and only if $d_1(y)=d_0(z)$.
Thus, the event ``(a) or (b)'' happens if and only if
$d(y)=d(z)$.

Now, suppose instead that $\alpha\ne0$ and that $\alpha\in E_\mathrm{span}^0$.
Then, (a) happens if and only if $d_1(y)-d_1(z)$ is even,
and (b) happens if and only if $d_1(y)-d_0(z)$ is even.
Thus, the event ``(a) or (b)'' happens if and only if
$d(y)-d(z)$ is even or $\#V$ is odd.
\end{proof}

We are finally ready to prove our third and last main result, which is a dual
analogue to \cref{th:seven}.
\maindualmodtwolinegraph
\begin{proof}
The first part of the proof is identical to first part of the proof of
\cref{th:maindualnormalfield}.

Let $\theta_\beta$ be the affine mapping from $X$ to $X^\ast$ defined by
$\theta_\beta(x) = \omega(x) + \beta$.
Note that $\theta_\beta^{-1}(\beta)=\ker\omega$, which is nonempty.
By \cref{lm:pressingstolamps}, $\beta$ and $\alpha$ belong to the same
$\transgrsimple{S}^\ast$-orbit only if $\theta_\beta^{-1}(\gamma)$ is nonempty too,
that is, if $\beta-\alpha$ belongs to $\im\omega$.

Suppose $\gamma-\beta\in\im\omega$. By \cref{lm:directsum}, $X^\ast=\im\delta\oplus E_\mathrm{span}^0$.
Let $\alpha$ be the part of $\beta$ that belongs to $E_\mathrm{span}^0$ in this direct sum. Then, $\delta^{-1}(\beta-\alpha)$
is nonempty.
By \cref{pr:deltapartialomega}, $\delta\circ\partial=\omega$, so $\gamma-\beta\in\im\delta$ and it follows that $\delta^{-1}(\gamma-\alpha)$ is nonempty too.
By \cref{lm:pressingstolampslinegraph},
$\beta$ and $\gamma$ belong to the same
$\transgrsimple{S}^\ast$-orbit if and only if $\delta^{-1}(\beta-\alpha)$ and
$\delta^{-1}(\gamma-\alpha)$ belong to the same
$\restrict{\transgr{\alpha}{S}}{X/\ker\omega}$-orbit.

First suppose $\beta\not\in\im\delta$. Then $\alpha\ne0$, and
by \cref{lm:sameorbitfromd}, $\delta^{-1}(\beta-\alpha)$ and
$\delta^{-1}(\gamma-\alpha)$ belong to the same orbit
if and only if $d(\delta^{-1}(\beta-\alpha))-d(\delta^{-1}(\gamma-\alpha))$
is even or $\#V$ is odd. If $\#V$ is even,
$d_0(y)$ has the same parity as $d_1(y)$ for any $y\in\generate{V}$, so
$d$ modulo 2 is a linear mapping from $\generate{V}$ to $K$.
Thus,
\[
d(\delta^{-1}(\beta-\alpha))-d(\delta^{-1}(\gamma-\alpha))
\equiv_2 d(\delta^{-1}(\beta-\gamma)),
\]
and we claim that this is zero modulo 2.
Since $\beta-\gamma$ belongs to $\im\omega$ and $\omega=\delta\circ\partial$,
there is an $x\in X$ such that
$\delta(\partial(x))=\beta-\gamma$, which implies that
$d(\delta^{-1}(\beta-\gamma))=d(\partial(x))$. Now, since $d(\partial(e))\equiv_2 0$ 
for any $e\in E$, $d(\partial(x))\equiv_2 0$ too. The claim is proven, and
we conclude that $\beta$ and $\gamma$ belong to the same orbit.

Now suppose instead that $\beta\in\im\delta$. Then $\alpha=0$, and
by \cref{lm:sameorbitfromd}, $\delta^{-1}(\beta-\alpha)=\delta^{-1}(\beta)$ and $\delta^{-1}(\gamma-\alpha)=\delta^{-1}(\gamma)$ belong to the same orbit if and only if they have the same $d$-value.
\end{proof}

\section{Handling multiple components in the ordinary case}\label{sec:nondual}
The three theorems listed in the introduction, \cref{th:brownhumphriesnormalfield,th:brownhumphriesFtwo,th:seven} all assume
that $G(S)$ is connected, and so do our dual analogues, \cref{th:maindualnormalfield,th:maindualmodtwoorthogonal,th:maindualmodtwolinegraph}.
For the dual case this is not a restriction, since
\cref{th:components} tells us how to handle multiple components, but what can we
say about the orbits of $\transgrsimple{S}$ acting on $X$ if $G(S)$ is not connected?

In this section we answer that question provided that $K\ne\bbF_2$.
\begin{theo}
Suppose $K\ne\bbF_2$, and let $S$ be a spanning subset of $X$.
Then, two elements
$x,y\in X$
belong to the same orbit of $\transgrsimple{S}$ if and only if
$y-x$ belongs to the span of the union of all $S_i$ such that
neither $\omega(x,S_i)$ nor $\omega(y,S_i)$ is $\{0\}$.
\end{theo}
\begin{proof}
Since $S$ spans $X$, we can write $x=\sum_{i\in I}x_i$, where $x_i$
belongs to $\Span S_i$. Let $J$ be the subset of $i\in I$ such that
neither $\omega(x,S_i)$ nor $\omega(y,S_i)$ is $\{0\}$.

First suppose $g(x)=y$ for some $g\in\transgrsimple{S}$.
We have
\[
g(x)=g({\textstyle\sum_{i\in I}x_i})=\sum_{i\in I}g(x_i).
\]
Let $I_x$ be the set of $i\in I$ such that $\omega(x,S_i)=\{0\}$
and let $I_y$ be the set of $i\in I$ such that $\omega(S_i,y)=\{0\}$.
For any $i\in I_x$ we have $g(x_i)=x_i$,
so $y-x=g(x)-x$ belongs to $\Span\bigcup_{i\in I\setminus I_x}S_i$.
By an analogous argument, $y-x=y-g^{-1}(y)$ belongs to
$\Span\bigcup_{i\in I\setminus I_y}S_i$, and we conclude that
$y-x$ belongs to the span of the union of all $S_j$ with $j\in J$.

For the converse, suppose $y-x$ belongs to the span of the
union of all
$S_j$ with $j\in J$. Then $y$ can be written as
$y=\sum_{i\in I}y_i$, where $y_i\in\Span S_i$ for any $i$
and $y_i=x_i$ for any $i\in I\setminus J$.
By the definition of $J$, for any $j\in J$
there is an $s\in S_j$ with $\omega(x,s)\ne 0$,
which implies that $\omega(x_j,s)\ne 0$, so $x_j$ does not
belong to the kernel of $\restrict{\omega}{\Span S_j}$,
the restriction of the form to $\Span S_j$. By the same argument,
$y_j$ does not belong to $\ker\restrict{\omega}{\Span S_j}$ either.
By \cref{th:brownhumphriesnormalfield}, $x_j$ and $y_j$ belong to the same
orbit of $\restrict{\transgrsimple{S_j}}{\Span(S_j)}$. Thus, there are $g_i\in\restrict{\transgrsimple{S_i}}{\Span S_i}$ such that $g_i(x_i)=y_i$ for any $i\in I$,
and by the group isomorphism result in \cref{th:components},
there is a unique $g\in\transgrsimple{S}$ such that
$\restrict{g}{\Span S_i}=g_i$ for any $i$. We have
\[
g(x)=g({\textstyle\sum_{i\in I}x_i})=\sum_{i\in I}g(x_i)=
\sum_{i\in I}g_i(x_i)=\sum_{i\in I}y_i=y.
\]
\end{proof}

\bibliographystyle{abbrv}
\bibliography{bib-file/numbers_game_mod_two}

\end{document}